\documentclass[11pt]{article}

\newcommand{\pr}[1]{\hbox{{\bf P}${}^{#1}$}}

\newtheorem{thm}{Theorem}[section]

\newtheorem{lem}[thm]{Lemma}
\newtheorem{cor}[thm]{Corollary}
\newtheorem{rem}[thm]{Remark}
\newtheorem{example}[thm]{Example}

\newcommand{\skipit}[1]{{}}
\newcommand{\prfend}{\hbox to7pt{\hfil}
\par\vskip-\baselineskip\hbox to\hsize
{\hfil\vbox {\hrule width6pt height6pt}}\vskip\baselineskip}
\newcommand{\putit}[4]{\noindent
\vbox to0in{\vskip#2in\hskip#1in\vbox{\hsize#4in #3}\vss}
\vskip-\baselineskip}
\newcommand{\theform}[3]{\noindent\vbox to0in{\vskip#2in
\hskip#1in\includegraphics{#3}
\vss}}
\newcommand{\tmprow}[2]{\hbox to\hsize{\hskip.3in\hbox to5.5in{\hbox to3in{#1\hfil}#2\hfil}\hss}\par}
\newcommand{\C}[1]{\hbox{$\cal #1$}}
\newcommand{\mytmprow}[4]{\hbox to\hsize{\hskip.2in\hbox to5.25in{\hbox to4.15in{\hbox to1.75in{\hbox to .6in{#1\hfil}#2\hfil}#3\hfil}#4\hfil}\hss}\par}
\newcommand{\myrow}[4]{\hbox to\hsize{\hskip.2in\hbox to5.25in{\hbox to3.35in{\hbox to1.75in{\hbox to .6in{#1\hfil}#2\hfil}#3\hfil}#4\hfil}\hss}\par}

\begin{document}
\title{Resolutions of ideals of any six fat points in \pr2}

\author{Elena Guardo\\
Dipartimento di Matematica e Informatica\\
Viale A. Doria 6, 95100\\
Catania\\
Italy\\
email: guardo@dmi.unict.it\\
WEB: http://www.dmi.unict.it/$\scriptstyle\sim$guardo
\and
Brian Harbourne\\
Department of Mathematics\\
University of Nebraska-Lincoln\\
Lincoln, NE 68588-0130\\
USA\\
email: bharbour@math.unl.edu\\
WEB: http://www.math.unl.edu/$\scriptstyle\sim$bharbour/}

\maketitle

\begin{abstract}
The graded Betti numbers of the minimal free resolution
(and also therefore the Hilbert function)
of the ideal of a fat point subscheme $Z$ of \pr2 are determined 
whenever $Z$ is supported at any 6 or fewer distinct points. 
All results hold over an
algebraically closed field $k$ of arbitrary characteristic.
\end{abstract}

\thanks{Acknowledgments: We would like to thank the GNSAGA, the NSA, the NSF and the 
Department of Mathematics at UNL for their support of the authors' research and
of E. Guardo's visits in 2003, 2004 and 2005 while this work was carried out.
We would also like to thank Jeremy Martin, Hal Schenck and Hannah Markwig
for helpful comments.}

\section{Introduction}\label{intro}

We begin by describing the problem we solve here, using terminology
familiar to experts. Those readers not already familiar
with the jargon can rest easy, since we will recall 
what the terms mean in section \ref{bkgd}. 

Given general points $p_1,\ldots,p_n$ of \pr2 and arbitrary positive
integers $m_i$, it is an open problem to determine the
graded Betti numbers for the minimal free resolution
of the ideal $I(Z)$ of the fat point subscheme $Z=m_1p_1+\cdots+m_np_n$.
It is even an open problem to determine just the Hilbert function of $I(Z)$.
Partly because of the difficulty of these problems and because
a standard approach to them involves considering special configurations
of points, and partly because of the intrinsic interest, there 
has been growing interest in these problems 
not only for general points but also when the points need not be general,
both in the plane and in higher dimensions
(see, for example, 
\cite{refBGVTone},
\cite{refBGVTtwo},
\cite{refCat},
\cite{refFHL},
\cite{reffranc},
\cite{refGMS},
\cite{refGVT},
\cite{refGVTb},
\cite{refTAMS},
\cite{refBHProc},
\cite{refARS},
\cite{refFreeRes},
\cite{refHR}).

In particular, \cite{refGMS} raises the question of finding
all Hilbert functions and graded Betti numbers for ideals of double point
subschemes of the plane; i.e., for $2p_1+\cdots+2p_n\subset\pr2$,
for all possible configurations of the points $p_i$. 
As \cite{refGMS} discusses, the Hilbert functions which occur for simple point
subschemes $Z=p_1+\cdots+p_n\subset\pr2$ are known
for all possible configurations of the points $p_i$; 
one goal that \cite{refGMS} works toward is
to find all Hilbert functions occurring for double point subschemes
$2Z=2p_1+\cdots+2p_n\subset\pr2$ such that the support scheme
$Z=p_1+\cdots+p_n\subset\pr2$ has given Hilbert function.
While \cite{refGMS} shows that for each Hilbert function of simple
points there is a Hilbert function which
in each degree has minimal value, it leaves unsolved the problem
of how to actually find this minimal Hilbert function, even for
small values of $n$ (such as $n=6$), and it raises
the question of whether there is also a maximal Hilbert function.
(It is worth mentioning that while we talk about the Hilbert function
of the ideal $I(Z)$, \cite{refGMS} talks about the Hilbert function
of the quotient ring $R/I(Z)$, where $R$ is the homogeneous
coordinate ring of \pr2. Thus what is for us a maximal Hilbert function
is for \cite{refGMS} a minimal Hilbert function.)

We answer all of these questions for the case of 6 points of \pr2 (see 
section \ref{distpnts}).
Moreso, we give a general approach for answering
any problems of the kinds raised in \cite{refGMS}, for {\it any\/} fat point
subschemes of \pr2 with support at 6 points, regardless of the
multiplicities $m_i$. More precisely, 
define a configuration type of $n$ points by 
requiring that sets $\{p_1,\ldots,p_n\}\subset\pr2$
and $\{p'_1,\ldots,p'_n\}\subset\pr2$ of distinct points 
have the same {\it configuration type\/} 
if and only if, after reordering the points $p'_i$ if need be,
the ideals of $Z=m_1p_1+\cdots+m_np_n$
and $Z'=m_1p'_1+\cdots+m_np'_n$ have the same Hilbert function
for every choice of the nonnegative integers $m_i$.
We show not only that the set of all
configurations of 6 points of \pr2 fall into only 11 different
types (see Corollary \ref{typesvsneg} and 
section \ref{distpnts}), but that ideals of
any two subschemes $Z=m_1p_1+\cdots+m_6p_6$
and $Z'=m_1p'_1+\cdots+m_6p'_6$
whose points have the same type 
also have the same graded Betti numbers (see 
Theorem \ref{distpntsThm} and Example \ref{distptsex}).
Our method also allows us to write down the Hilbert function
and graded Betti numbers for any $Z=m_1p_1+\cdots+m_6p_6$,
given only the coefficients $m_i$ and given the configuration type
of the points with respect to a specific ordering of the points. 
(Figure 1 shows the 11 different configuration types 
of 6 points. Thus type 1 consists of 6 general points; for
type 2, three of the points are collinear, etc. Type
11 has all six points on an irreducible conic.)

What is new here is the explicit enumeration of the 11 types
(this is easy), and the determination of the graded Betti numbers 
(this is where most of the effort of this paper lies).
It follows from our main result, Theorem \ref{distpntsThm}, 
that numerical Bezout considerations (as discussed in 
Remark \ref{HilbFuncRem} and demonstrated in Example \ref{distptsex})
suffice to determine the graded Betti numbers 
of a fat point ideal supported at any 6 distinct points 
of \pr2. (By numerical Bezout considerations we are referring to
the version of Bezout's theorem that tells us that two effective divisors
$C$ and $D$ on an algebraic surface must have a common component
if their intersection $C\cdot D$ is negative. We give a procedure
for computing the graded Betti numbers that depends
only on computing intersections of divisor classes on
a blow up of \pr2, which amounts to taking dot products
of integer entry vectors. This procedure is easy to carry 
out by hand, as shown by Example \ref{distptsex}. An awk script implementing it can 
be run over the web by visiting\par
\noindent{\tt http://www.math.unl.edu/$\scriptstyle\sim$bharbour/6ptres/6reswebsite.html} .)

The facts that, for any $n\le 8$ points of the plane, numerical Bezout considerations determine
the Hilbert function of any fat point subscheme supported at those points, 
and that there are only finitely many
different configuration types of $n\le 8$ points, follow from the main result
of \cite{refBHProc}. However, these facts seem not to be widely
recognized (the authors of \cite{refGMS}, for example, were
quite interested when we mentioned this to them),
perhaps because the finite set of configurations has never
been explicitly written down. 

Enumerating these finitely many types for $7\le n\le8$
takes considerably more effort than doing so for $n=6$;
in a not yet written preprint, Geramita, Harbourne and Migliore
find 29 types of distinct points for $n=7$ and 146 for $n=8$.
Determining how the graded Betti numbers
behave will be much more difficult, both because of the many
cases that need to be considered, and because the behavior of the graded
Betti numbers is more subtle (see  \cite{refCJM} and \cite{refFHH},
which work out the graded Betti numbers for 7 and 8 general points respectively,
versus \cite{refFi}, which works out the case of 6 general points).
Moreover, for $n>8$ points, the number of types
is infinite. (For example, just by taking points in
various configurations on a smooth non-supersingular
plane cubic curve, for any positive integer $r$ one can 
by Proposition 1.2 of
\cite{refTAMS} arrange for the Hilbert function
of $I(mp_1+\cdots+mp_9)$ in degree $t=3m$ to 
be $\lfloor m/r\rfloor+1$. Thus the number of Hilbert functions
increases with $m$, so for $m$ large enough no given
finite set of types will be sufficient to encompass all of them.)

We now briefly discuss additional background for our
work in this paper. 
The Hilbert function for ideals $I(Z)$ of fat point subschemes
$Z\subset \pr2$ supported at $n\le 9$ general points 
is well known; see, for example, Nagata \cite{refNtwo},
or, for $n=6$, Giuffrida \cite{refGf}. For $n>9$ general points, 
the problem of finding the Hilbert function of $I(Z)$ has been solved
only in special cases. As mentioned above, the problem of finding 
the Hilbert function of $I(Z)$ as long as $Z$ has support
at $n\le8$ points, even possibly infinitely near, 
was solved, in principle, in \cite{refBHProc}, without
however classifying the possible configuration types.

\vskip3.05in


\putit{1.6}{-2.7}{Type 1 \hfil Type 2\hfil Type 3}{3}

\putit{1.6}{-2}{Type 4 \hfil Type 5\hfil Type 6}{3}

\putit{1.6}{-1.1}{Type 7 \hfil Type 8\hfil Type 9}{3}

\putit{1.8}{-.15}{Type 10 \hfil Type 11}{3}

\theform{.7}{2.13}{"11configsNoText.ps" hscale=55 vscale=55}

\hbox to\hsize{\hfil Figure 1\hfil}  
\vskip.1in

A logical next step is to determine the graded Betti numbers for
minimal free resolutions of ideals of fat point subschemes
in \pr2 supported at any configuration of points. Previous results have been given
in various cases. The first results are due to Catalisano \cite{refCat},
who determined the minimal free resolutions for fat point subschemes
supported at distinct points on an irreducible plane conic. 
The case that the conic is not irreducible or the points are possibly infinitely near
was handled in \cite{refFreeRes}.
(Since a connected curve of degree at most 2 in any projective space
lies in a plane in that space, by applying \cite{refFHL} the results of
\cite{refCat} and \cite{refFreeRes} actually also give
the Hilbert function and graded Betti numbers
for fat points in projective space of any dimension,
as long as the support of the points is contained
in a connected curve of degree at most 2.)
Various cases in which the points of the support are contained
in complete intersections in  
\pr2 are studied in \cite{refBGVTone}, \cite{refBGVTtwo}, \cite{refGVT} and  \cite{refGVTb}.
Additional special configurations are handled in \cite{refGMS}, but only
in case of points of multiplicity 2.

Since any five points lie on a smooth conic,
Catalisano's result handles the case of fat point 
subschemes supported at five general points.
The case of 6 general points was worked out by Fitchett \cite{refFi}.
For the case of seven general points, see \cite{refCJM}, and for eight general
points, see \cite{refFHH}. Numerous special cases for 9 or more general points
have been done (for $n\ge9$ general points of multiplicity 1, 2 or 3,
see \cite{refGGR}, \cite{refId} and \cite{refGI}, respectively; for $n$ general points
of multiplicity $m$ when $m$ is not too small and $n$ is an even square,
in light of \cite{refE}, see \cite{refHHF}; additional cases are handled by
\cite{refHR}). The problem for general points is otherwise open. There is a conjecture for
the Hilbert function of the ideal of any fat point subscheme of \pr2 supported
at general points (see \cite{refSurvey} for a discussion), and 
there are conjectures in special cases for resolutions 
(see \cite{refIGC} and \cite{refHHF}), but so far no
general conjecture for the resolution has been posed.

In this paper, we extend \cite{refCat} and \cite{refFreeRes}
to the case of any 6 distinct points of \pr2. 
Our approach involves a case by case analysis for the different configuration
types of 6 points in \pr2,
depending on finding sets of generators of the cone of nef divisor
classes on the surface $X$ obtained by blowing up
the 6 points. At first glance
verifying our result even for a single configuration of points
would seem to require checking an infinite number of cases,
since there are infinitely many nef divisor classes.
The fact that our methods make the problem tractable is of
interest in its own right. 

\section{Background}\label{bkgd}

We begin by discussing our methods in more detail.
So let $p_1,\ldots,p_n$ be distinct points
of \pr2. Given nonnegative integers $m_i$, the fat point
subscheme $Z=m_1p_1+\cdots+m_np_n\subset \pr2$ is, by definition, defined by the ideal
$I(Z)=I(p_1)^{m_1}\cap\cdots\cap I(p_n)^{m_n}$, where
$I(p_i)\subset R=k[\pr2]$ is the ideal generated by all forms
(in the polynomial ring $R$ in three variables over the field $k$)
vanishing at $p_i$. The {\it support\/} of $Z$ consists of 
the points $p_i$ for which $m_i$ is positive.

The minimal free resolution of $I(Z)$ is an exact sequence of the form
$$0\to F_1\to F_0\to I(Z)\to 0$$
where each $F_i$ is a free graded $R$-module, where the grading is with respect to the usual
grading of $R$ by degree, and all entries of the matrix defining
the homomorphism $F_1\to F_0$ are homogeneous polynomials in $R$ of degree
at least 1. To determine $F_0$ up to 
graded isomorphism, it is enough to determine the dimensions
of the cokernels of the multiplication
maps $\mu_{Z,i} : I(Z)_i\otimes R_1\to I(Z)_{i+1}$
for each $i\ge0$, where, given a graded $R$-module $M$, $M_t$
denotes the graded component of degree $t$. If we denote $\hbox{dim cok}(\mu_{Z,i-1})$
by $t_i$, then $F_0=\oplus_{i>0} R[-i]^{t_i}$, where $R[-i]$ is the free
graded $R$-module of rank 1 with a shift in degrees given by $R[-i]_j=R_{j-i}$.
The Hilbert functions of $I(Z)$ and $F_0$ then determine $F_1$ up 
to graded isomorphism. In fact, if we denote the Hilbert function
of $Z$ by $h_Z$ (i.e., $h_Z(i)=\hbox{dim }I(Z)_i$), and if $\Delta$ denotes the 
difference operator (i.e., $\Delta h_Z(i)=h_Z(i)-h_Z(i-1)$),
then $F_1=\oplus_{i>0}R[-i]^{s_i}$, where $s_i=t_i-(\Delta^3h_Z)(i)$ (see \cite{refFHH},
p. 685).

Thus to determine $F_0$ and $F_1$ it is enough to determine the Hilbert
function of $I(Z)$ and the rank of $\mu_{Z,i}$ for each $i$. 
The Hilbert function of $I(Z)$ can be obtained by applying the result of \cite{refBHProc}.
It follows from Theorem \ref{distpntsThm} that the ranks of the $\mu$ can be found
by a maximal rank criterion, as we now explain. 

Given $Z$, let $\alpha(Z)$ be the least degree $j$ such that $h_Z(j)>0$;
i.e., such that $I(Z)_j\ne0$. For each $t\ge \alpha(Z)$, 
let $\gamma(Z,t)$ be the gcd of $I(Z)_t$. Thus 
$\gamma(Z,t)$ is a homogeneous form of some degree $d_{Z,t}$. If $d_{Z,t}=0$,
it is convenient to set $\gamma(Z,t)=1$, but if $d_{Z,t}>0$, then
$\gamma(Z,t)$ defines a plane curve $C=C_{Z,t}$ of degree $d_{Z,t}$. Let $m'_i$ be 
the multiplicity $\hbox{mult}_{p_i}(C)$ of the curve at the point $p_i$.
Thus we get a fat points subscheme $Z_t^-=m'_1p_1+\cdots+m'_np_n$.
Let $Z_t^+=(m_1-m'_1)_+p_1+\cdots+(m_n-m'_n)_+p_n$, where for 
any integer $m$, $m_+=\hbox{max }(0,m)$. Then clearly
$I(Z)_t=\gamma(Z,t)I(Z_t^+)_{t-d_{Z,t}}$. 

For $n\le 8$ and $t\ge \alpha(Z)$, it is known that 
$$\hbox{dim}(I(Z_t^+)_{t-d_{Z,t}})=
{t-d_{Z,t}+2\choose 2}-\sum_i{(m_i-m'_i)_++1\choose 2},$$
as a consequence of the fact that a nef divisor $F$ on a blow up $X$ of \pr2 at $n\le 8$
points has $h^1(X, \C O_X(F))=0=h^2(X, \C O_X(F))$
\cite{refBHProc}. For $n\le 8$, as we discuss in more detail below in 
Remark \ref{HilbFuncRem},
one can determine $\alpha(Z)$ using purely numerical Bezout considerations,
and for each $t\ge\alpha(Z)$, one can also determine $Z_t^-=m'_1p_1+\cdots+m'_np_n$
and $d_{Z,t}$ purely numerically, from Bezout considerations.
(In order to determine these quantities
in the case of $n=6$ distinct points, in addition to having the coefficients $m_i$,
one needs to know only the configuration type with respect to
a specific ordering of the points; i.e., one needs to
know only whenever there is a line going through three or 
more of the points $p_i$, and which 
points those are, and if there is a conic going through all 6 points.)

Given that we can determine the Hilbert function of the ideal $I(Z)$,
to determine the graded Betti numbers $t_i$ and $s_i$ of the resolution, therefore,
it is enough to determine $t_i$ for each $i$. 
Since we know the Hilbert function, we know $\alpha(Z)$ and
clearly, $t_i=0$ for $i<\alpha(Z)$, and $t_i=h_Z(\alpha(Z))$ for $i=\alpha(Z)$.
If $i$ is large enough, the Hilbert function and Hilbert polynomial coincide;
i.e., we will have
$\hbox{dim}(I(Z)_i)={i+2\choose 2}-\sum_i{m_i+1\choose 2}$. Let
$\tau(Z)$ be the least $i$ such that this holds, and let $\sigma(Z)=\tau(Z)+1$.
Regularity considerations \cite{refDGM} then imply that 
$t_i=0$ for $i>\sigma(Z)$.

So assume $\alpha(Z) \le i <\sigma(Z)$. 
Since $I(Z)_i=\gamma(Z,i)I(Z_i^+)_{i-d_{Z,i}}$ for $i\ge\alpha(Z)$,
multiplying by $\gamma(Z,i)$ gives an inclusion $I(Z_i^+)_{i-d_{Z,i}+1}\subset I(Z)_{i+1}$
and a vector space isomorphism
between the images of $\mu_{Z_i^+,i-d_{Z,i}}$ and $\mu_{Z,i}$.
From the inclusions $\hbox{Im}(\mu_{Z_i^+,i-d_{Z,i}})\subset I(Z_i^+)_{i-d_{Z,i}+1}\subset I(Z)_{i+1}$
it now follows that
 
$$t_{i+1}=\hbox{dim cok}(\mu_{Z,i}) = \hbox{dim cok}(\mu_{Z_i^+,i-d_{Z,i}})
+ (h_Z(i+1) - h_{Z_i^+}(i-d_{Z,i}+1)).$$

Since $\gamma(Z_i^+,i-d_{Z,i})=1$, and assuming that we can determine Hilbert functions,
this reduces the problem of computing
$\hbox{dim cok}(\mu_{Z,i})$ for an arbitrary $Z$ in degrees $i\ge \alpha(Z)$ to the problem of computing
$\hbox{dim cok}(\mu_{Z,i})$ for an arbitrary $Z$ but only in degrees $i\ge\alpha(Z)$ such that
$\gamma(Z,i)=1$.
This is what we do. Our main result, Theorem \ref{distpntsThm}, essentially says that if
$Z$ has support at any 6 distinct points of \pr2, and if 
$i\ge \alpha(Z)$ is such that $\gamma(Z,i)=1$, 
then $\mu_{Z,i}$ has maximal rank (meaning that $\mu_{Z,i}$
is either injective or surjective and hence 
$t_{i+1}$ is either $h_Z(i+1)-3h_Z(i)$ or 0, respectively).
Since $\gamma(Z_i^+,i-d_{Z,i})=1$, it follows that
$\mu_{Z_i^+,i-d_{Z,i}}$ has maximal rank, and hence 
everything on the right hand side of the displayed formula
above is in terms of Hilbert functions of fat points supported
at the given 6 points. Computing those Hilbert functions
thus computes $\hbox{dim cok}(\mu_{Z,i})$. 

In order to compute the graded Betti numbers for the minimal
free resolution of fat point subschemes $Z$ with support
at 6 points, we thus need to determine their Hilbert functions
and, for each degree $i$, we need to determine
$Z_i^+$ and the degree of $\gamma(Z,i)$. The easiest
context in which this can be done involves the intersection
theory on the surface obtained by blowing up the points.
This will also be the context we use to study the rank of
$\mu_{Z,i}$.

Let $\pi:X\to \pr2$ be the birational morphism obtained by blowing up
distinct points $p_1,\ldots,p_n$ of \pr2. 
Let $\hbox{Cl}(X)$ be the divisor class group of $X$.
Let $E_0$ be the pullback to $X$ of the class of a line on \pr2,
and let $E_1,\ldots,E_n$ be the classes of the exceptional 
divisors of the blow ups of $p_1,\ldots,p_n$. 
Then $\hbox{Cl}(X)$ is formally just a free abelian
group with a preferred orthogonal basis $E_0,\ldots,E_n$.
This basis is called an {\it exceptional configuration}.
(The bilinear form on $\hbox{Cl}(X)$ is given by $E_i\cdot E_j=0$
for all $i\ne j$, $E_0^2=1$ and $E_i^2=-1$ for $i>0$.)
We are mainly interested in the case that $n=6$; hereafter,
we will often but not always assume that $n=6$.

Problems involving fat points with support at points
$p_1,\ldots,p_n$ on \pr2 can be translated
to problems involving divisors on $X$.
Given $Z$ and $t$, the vector space $I(Z)_t$ 
is a vector subspace of the space of sections
$H^0(\pr2,\C O_{\pr2}(t))$. The latter is referred
to as a complete linear system; $I(Z)_t$ is typically
a proper subspace, in which case it is referred to as 
an incomplete linear system. However, we can
associate to $Z=m_1p_1+\cdots+m_np_n$ and $t$ the 
divisor class $F(Z,t)=tE_0-m_1E_1-\cdots-m_nE_n$ on $X$,
in which case $I(Z)_t$ can be canonically identified (as a vector space) with
the complete linear system $H^0(X, \C O_X(F(Z,t)))$.

Given a divisor or divisor class $F$ on $X$, it will be
convenient to write $h^i(X, F)$ in place of 
$h^i(X, \C O_X(F))$, and we will refer to a divisor class 
$F$ as {\it effective\/} if $h^0(X, F)>0$; i.e., if it is the class of
an effective divisor.
In particular, $\hbox{dim } I(Z)_t = h^0(X, F(Z,t))$
for all $Z$ and $t$, and the ranks of 
$\mu_{Z,t}$ and $\mu_{F(Z,t)}$ are equal,
where
$$\mu_{F(Z,t)}:H^0(X, F(Z,t))\otimes
H^0(X,  E_0)\to H^0(X, F(Z,t)+E_0)$$
is the natural map given by multiplication.

Whenever $N$ is a
prime divisor (i.e., a reduced irreducible curve)
such that $F(Z,t)\cdot N<0$, we have
$h^0(X, F(Z,t))=h^0(X, M)$, where $M=F(Z,t)-N$. Moreover, clearly
the kernels of $\mu_{F(Z,t)}$ and $\mu_M$ have the same dimension,
so if we can compute $h^0$ for arbitrary divisors on $X$, finding the rank
of $\mu_{F(Z,t)}$ is equivalent to doing so for $\mu_M$.
If we have a complete list of prime divisors $N$ of negative
self-intersection, then whenever $F(Z,t)$ is effective,
we can subtract off prime divisors of negative self-intersection
to obtain an effective class $M$ which is {\it nef}
(meaning that $M\cdot D\ge0$ for every effective divisor $D$),
in which case 
$h^0(X, F(Z,t))=h^0(X, M)$ and
the kernels of $\mu_{F(Z,t)}$ and $\mu_M$ have the
same dimension, thereby reducing the problem to the case of computing
$h^0(X, M)$ and ranks of $\mu_M$ only when $M$ is nef.

This is very helpful, since for $n\le 8$, $h^1(X, M)=0=h^2(X, M)$
whenever $M$ is nef (\cite{refBHProc})
and hence $h^0(X, M)=(M^2-K_X\cdot M)/2+1$ by Riemann-Roch.
Thus for $n\le 8$, the Hilbert function of $I(m_1p_1+\cdots+m_np_n)$ is completely 
determined by the coefficients $m_i$ and by the set of classes of prime divisors of negative 
self-intersection on the surface $X$ obtained by blowing up the points $p_i$. 
(For $n\ge9$, this is no longer true.
This is because $h^1(X,\C O_X(F))=0$ can fail for nef divisors
when $n\ge9$, as shown by considering a general pencil of cubics.)

But whereas $\mu_M$ is always surjective
for nef divisors $M$ for any $n\le5$ distinct (or
even possibly infinitely near)
points \cite{refFreeRes},
$\mu_M$ can fail to have maximal rank for nef divisors
when $n\ge7$ \cite{refCJM}, even for $n$ general points. 
However, for $n=6$ general points,
$\mu_M$ always at least has maximal rank when $M$ is nef
\cite{refFi}. This leaves open the question of whether
$\mu_M$ may fail to have maximal rank for some nef $M$ for 
some particular choice of $n=6$
distinct points; we show
that $\mu_M$ has maximal rank for any nef $M$
for all choices of the points $p_i$.

We begin by determining the subset $\hbox{NEG}(X)\subset\hbox{Cl}(X)$
of divisor classes of effective reduced irreducible
divisors of negative self-intersection. Among all 6 point blow ups $X$
of \pr2, it turns out there are 
only finitely many possible subsets $\hbox{NEG}(X)$,
and $\hbox{NEG}(X)$ is itself always finite. (By Corollary \ref{typesvsneg},
up to reordering the points, the possible subsets
$\hbox{NEG}(X)$ correspond bijectively with the configuration types
of Figure 1.)  
We can then obtain our result by an analysis for each possible 
subset $\hbox{NEG}(X)$. As a practical matter, it is easier
to consider the subset 
$$\hbox{neg}(X)=\{C\in \hbox{NEG}(X)\hbox{ : }C^2<-1\},$$ 
since $\hbox{neg}(X)$ is a proper  
(and usually substantially smaller) subset
of $\hbox{NEG}(X)$, but $\hbox{neg}(X)$
determines $\hbox{NEG}(X)$, by Remark \ref{minustwoRem}. In fact,
the elements of $\hbox{neg}(X)$ correspond to the 
curves displayed in Figure 1. (For example, 
for configuration type 1, $\hbox{neg}(X)$ is empty,
for configuration type 2, $\hbox{neg}(X)$ consists of the
divisor class of the proper transform of the line through the three collinear
points, etc.)

While $\hbox{NEG}(X)$ and $\hbox{neg}(X)$ depend
on the particular points $p_i$, we now define a fixed finite
subset of the divisor class group $\hbox{Cl}(X)$ which contains them.
Consider $\C B\cup\C L\cup\C Q$, where $\C B=\{E_i: i>0\}$
($\C B$ here is for {\it blow up} of a point),
$\C L=\{E_0-E_{i_1}-\cdots-E_{i_r}: r\ge 2, 0<i_1<\cdots<i_r\le 6\}$
($\C L$ here is for points on a {\it line}), and
$\C Q=\{2E_0-E_{i_1}-\cdots-E_{i_r}: r\ge 5, 0<i_1<\cdots<i_r\le 6\}$
($\C Q$ here is for points on a conic, defined by a {\it quadratic} equation).

The next result, which is well known but hard to cite in the form we need,
shows that there are only finitely many possibilities for
$\hbox{NEG}(X)$, since it is a subset of $\C B\cup\C L\cup\C Q$.
(The finiteness remains true as long as $n<9$ but can fail for 
$n\ge9$. In addition, more possibilities
occur than the ones listed here if $n$ is 7 or 8.)

\begin{lem}\label{NEGisfinite}
Let $X$ be obtained by blowing up
6 distinct points of \pr2. Then the following hold:
\begin{description}
\item[(a)] $\hbox{NEG}(X)\subset \C B\cup\C L\cup\C Q$,
and every class in $\hbox{NEG}(X)$ is the class 
of a smooth rational curve;
\item[(b)] for any nef $F\in\hbox{Cl}(X)$, $F$ is effective
(hence $h^2(X, F)=0$), $|F|$ is base point free,
$h^0(X, F) = (F^2-K_X\cdot F)/2 + 1$ and $h^1(X, F)=0$;
\item[(c)] $\hbox{NEG}(X)$ generates the subsemigroup 
$\hbox{EFF}(X)\subset\hbox{Cl}(X)$ of classes of effective divisors;
and
\item[(d)] any class $F$ is nef if and only if 
$F\cdot C\ge0$ for all $C\in \hbox{NEG}(X)$.
\end{description}
\end{lem}

\noindent {\em Proof}. Riemann-Roch for a smooth rational surface
$X$ states that $h^0(X, A)-h^1(X, A)+h^2(X,A) = 
(A^2-K_X\cdot A)/2 + 1$
holds for any divisor class $A$. 
Also, $-K_X=3E_0-E_1-\cdots-E_6$, so $-K_X\cdot E_0=3$.
If $F$ is effective, then $F\cdot E_0\ge0$, since
$E_0$ is nef. 
(The reason $E_0$ is nef is that it is the class of an irreducible divisor
of nonnegative self-intersection, hence any effective
divisor meets it nonnegatively. More generally, 
any effective divisor which meets each of its 
components nonnegatively is nef.) By duality,
$h^2(X, F)=h^0(X, K_X-F)$, and $h^0(X, K_X-F)=0$
since $-K_X\cdot E_0=3$, hence $(K_X-F)\cdot E_0<0$.
This verifies the parenthetical remark in part (b).
Similarly, $h^2(X,-K_X)=0$, so we have  $h^0(X, -K_X) = 
K_X^2 +1 + h^1(X, -K_X)$, but for us $K_X^2=3$, so
$h^0(X, -K_X)=4+h^1(X, -K_X)$. Thus
$-K_X$ is the class of an effective divisor, say $D$.
Moreover, the subgroup $K_X^\perp\subset \hbox{Cl}(X)$
of all classes orthogonal to $-K_X$ is negative definite.
This is easy to see since 
the classes $E_1-E_2$, $E_1+E_2-2E_3$, 
$E_1+E_2+E_3-3E_4$, $E_1+E_2+E_3+E_4-4E_5$ and
$2E_0-E_1-\cdots-E_6$ have negative self-intersection
but are linearly independent and pairwise orthogonal,
hence give an orthogonal basis
of $K_X^\perp$ over the rationals. On the other hand,
it is not hard to check that
$E_0-E_1-E_2-E_3$, $E_1-E_2$, $\dots$, $E_5-E_6$
give a ${\bf Z}$-basis for $K_X^\perp$, and since
each basis element has self-intersection $-2$, it follows
that $A^2$ is even for every $A\in K_X^\perp$. I.e.,
$K_X^\perp$ is even and negative definite.

To justify (a), let $C$ be the class of a reduced irreducible divisor
on $X$, with $C^2<0$. Since $E_0$ is nef, we know $C\cdot E_0\ge0$.
If $C\cdot E_0=0$, then $C$ must be a component of one of the 
$E_i$, hence $C\in \C B$, since each $E_i$ is 
reduced and irreducible. If $C\cdot E_0=1$, then
$C$ is the proper transform of a line in \pr2, so
$C\in\C L$. If $C\cdot E_0=2$, then $C$ is the proper
transform of a smooth conic in \pr2, so
$C\in \C Q$. By explicitly applying adjunction $C^2+C\cdot K_X=2g-2$, 
where $g$ is the ({\it a priori\/} arithmetic) genus of $C$, any
$C\in \C B\cup\C L\cup\C Q$ which is the class of a prime divisor
has $g=0$ and so is the class of a smooth rational curve.

Now it suffices to show that we cannot have
$C\cdot E_0>2$.
If $C\cdot D < 0$, then $C$ is the class of an irreducible component
of $D$, hence $E_0\cdot (-K_X-C)\ge0$, so $E_0\cdot C\le 3$.
If $C\cdot E_0=3$, then $C$ is the proper
transform of an irreducible plane cubic. But an irreducible
plane cubic has at most one singular point, which must be of
multiplicity 2. Thus its proper transform is either
$3E_0-E_{i_1}-\cdots-E_{i_r}$, with $0<i_1<\cdots<i_r\le6$,
or $3E_0-2E_{i_1}-E_{i_2}-\cdots-E_{i_r}$, with $0<i_2<\cdots<i_r\le6$
and $0<i_1<6$ such that $i_1\ne i_j$ for $j>1$. But in neither case
would we have $C^2<0$, so $C\cdot E_0\ge3$ cannot happen. 

Now say $C\cdot D\ge 0$. 
From adjunction, since $0\le C\cdot D = -K_X\cdot C$, we have
$-1\ge C^2\ge -2$, with $g=0$ in any case, hence $C$ is a smooth
rational curve. If $C\cdot D=0$, then $C\in K_X^\perp$ and
adjunction gives $C^2=-2$,
but since $K_X^\perp$ is negative definite, it has only finitely
many classes of self-intersection $-2$. One can show
that the only classes in $K_X^\perp$ of self-intersection $-2$
are $\pm E_i\pm E_j$, $0<i<j\le 6$, $\pm E_0\pm E_i\pm E_j\pm E_k$, 
$0<i<j<k\le 6$, and $\pm2E_0\pm E_1\pm\cdots\pm E_6$.
(To see this, assume that $A=aE_0-b_1E_1-\cdots-b_6E_6\in K_X^\perp$.
Thus $3a = b_1+\cdots+b_6$. Working over 
$\hbox{Cl}(X)\otimes_{\bf Z}{\bf Q}$, let $m=(b_1+\cdots+b_6)/6$,
so $a=2m$, and define $B=aE_0-m(E_1+\cdots+E_6)$. 
Then $B\cdot K_X=0$, but $A^2\le B^2=-2m^2$. If $A^2=-2$,
then we must have $a\le 2$, in order to have $m\le 1$.
Thus $a$ is either 0, 1 or 2, and now it is easy to enumerate solutions
$A^2=-2$.) Among these classes, only those in
$\C B\cup\C L\cup\C Q$ 
can be classes of prime divisors.
(This is because a prime divisor must, first, meet
$E_0$ nonnegatively, and second, 
when expressed as a linear combination $a_0E_0-\sum a_iE_i$,
if $a_j<0$ for some $j>0$, then it must be a component of 
$E_j$ and thus must be in $\C B$.)

If $C\cdot D > 0$, then 
$C^2=-1=K_X\cdot C$. Let $Y\to X$ be obtained by blowing up a
seventh, general point $p_7$. This morphism induces
an inclusion $\hbox{Cl}(X)\to\hbox{Cl}(Y)$. Then, arguing as above,
$K_Y^\perp$ is even and negative definite, and the only solutions
to $A^2=-2$ for $A\in K_Y^\perp$ are of the form
$\pm E_i\pm E_j$, $0<i<j\le 7$, $\pm E_0\pm E_i\pm E_j\pm E_k$, 
$0<i<j<k\le 7$, and $\pm2E_0\pm E_{i_1}\pm\cdots\pm E_{i_6}$,
$0<i_1<\cdots<i_6\le 7$.
Thus $C-E_7$ is in $K_Y^\perp$, with $(C-E_7)^2=-2$, since
$C\cdot E_7=0$, and (keeping in mind that $C$
is a prime divisor also on $Y$ and that $C\cdot K_X=-1$) it follows that
$C$ is among $E_i$, $0<i\le 6$, $E_0-E_i-E_j$, 
$0<i<j\le6$, and $2E_0-E_{i_1}-\cdots-E_{i_5}$,
$0<i_1<\cdots<i_5\le 6$. This finishes the proof of (a).

To prove (b), we have $h^1(X, F)=0$ and $h^2(X, F)=0$
by Theorem 8, \cite{refBHProc}. Thus 
$h^0(X, F) = (F^2-K_X\cdot F)/2 + 1$ follows by 
Riemann-Roch. But $F^2\ge0$ holds for nef divisors 
(Proposition 4, \cite{refBHProc}),
so $F$ is effective. To see that $|F|$ is base point free,
note that a nef divisor in $K_X^\perp$ must be 0.
Now apply Theorem III.1(a,b) of \cite{refARS} to see that
$|F|$ is fixed component free, and has a base point
only if $-K_X\cdot F = 1$, in which case, 
using $Y$ as above, we see that $F-E_7$ must be effective,
but $F-E_7\in K_Y^\perp$, so $F^2-1=(F-E_7)^2\le0$. But
$(F-E_7)^2=0$ implies $F-E_7=0$, which is impossible
since then $0=F\cdot E_7=E_7^2=-1$.
Thus $0>(F-E_7)^2=F^2-1$, so
$F^2=0$. However, we also have $-K_X\cdot F=1$,
which contradicts $h^0(X, F) = (F^2-K_X\cdot F)/2 + 1$,
since $h^0(X, F)$ must be an integer. Thus we cannot have
$-K_X\cdot F = 1$ if $F$ is nef.

Consider (c). Let $G$ be the class of an effective
divisor. We can write $G=N+F$, where $N$ is the fixed part of
$|G|$, and $F$ is nef. Note that no component
of $N$ can be nef, since nef divisors (in our situation)
are base point free, whereas components of $N$ are fixed. Thus the 
class of every component of $N$ is in $\hbox{NEG}(X)$. 
Now, if a class $F=a_0E_0-a_1E_1-\cdots
-a_6E_6$ is nef for a particular 
set of distinct points $p_i$,
then it remains nef when the points $p_i$ are general, and if $F$ is effective
when the points are general, it was effective to begin with.
(This is because by semicontinuity the effective subsemigroup can never get smaller
as the points are specialized, so the nef cone can never enlarge.) 
And if the points $p_i$ are general, then $\hbox{NEG}(X)$ consists
of the exceptional classes; i.e., the classes $E_i$, $i>0$,
$E_0-E_i-E_j$, $0<i<j\le6$, and $2E_0-E_{i_1}-\cdots-E_{i_5}$,
$0<i_1<\cdots<i_5\le 6$. It follows from \cite{refTAMS}, that
the class of every effective divisor is a nonnegative sum
of exceptional classes. (The results of \cite{refTAMS}
show that it is enough to show that
$E_0$, $E_0-E_1$, $2E_0-E_1-E_2$, and $3E_0-E_1-\cdots E_j$,
$3\le j\le 6$ are, but this is easy; for example, 
$E_0=(E_0-E_1-E_2)+E_1+E_2$.) 
Thus given a class $F$ which is nef for a given set of
points $p_i$, $F-E$ is effective
for some $E$ among the classes $E_i$, $i>0$,
$E_0-E_i-E_j$, $0<i<j\le6$, and $2E_0-E_{i_1}-\cdots-E_{i_5}$,
$0<i_1<\cdots<i_5\le 6$. If $E$ is a prime divisor, then 
$E\in \hbox{NEG}(X)$. If not, then $E\cdot N'<0$ for some 
$N'\in \hbox{NEG}(X)$ (otherwise, $E$ is nef, hence
$h^0(X, E)=(E^2-K_X\cdot E)/2+1=1$, but also
$|E|$ must be base point free, hence $h^0(X, E)>1$).

Thus either way there is an $N'\in \hbox{NEG}(X)$ such that
$F-N'$ is effective. 
By replacing $F$ by $F-N'$ and repeating the process,
we eventually reach the case that $F=0$, 
hence any effective divisor is a sum of 
elements of $\hbox{NEG}(X)$.

Finally, we prove (d). To show $F$ is nef, we just need to show
that $F\cdot C\ge 0$ for each class $C$ of an effective divisor.
But each such $C$ is a nonnegative sum of classes
in $\hbox{NEG}(X)$ and any class in $\hbox{NEG}(X)$ is 
the class of an effective divisor. It follows that
$F\cdot C\ge 0$ for the class $C$ of an effective divisor
if and only if $F\cdot C\ge 0$ for every $C\in \hbox{NEG}(X)$.
\prfend

\begin{rem}\label{minustwoRem}\rm We now show how
$\hbox{neg}(X)$ determines $\hbox{NEG}(X)$.
In fact, 
$$\hbox{NEG}(X)=\hbox{neg}(X)\cup
\{C\in\C B\cup\C L\cup\C Q\hbox{ $|$ }C^2=-1,
C\cdot D\ge0\ \forall\ D\in \hbox{neg}(X)\}.$$
The forward inclusion follows from 
Lemma \ref{NEGisfinite}(a).
For the reverse, say $C^2=-1$ for some
$C\in \C B\cup\C L\cup\C Q$. It is easy to check
case by case that each such $C$ is effective, hence 
$C\cdot C'<0$ for some $C'\in\hbox{NEG}(X)$.  
Given that $C\cdot D\ge0$ for all $D\in\hbox{neg}(X)$,
then $C'\in\hbox{NEG}(X)-\hbox{neg}(X)$.
But any two distinct elements of $\C B\cup\C L\cup\C Q$
of self-intersection $-1$ meet nonnegatively, hence
$C=C'\in\hbox{NEG}(X)$.
\end{rem}

By Lemma \ref{NEGisfinite} and Remark \ref{minustwoRem} it follows that
specifying $\hbox{neg}(X)$ as a subset of $\C L\cup\C Q$
is equivalent to specifying the configuration type of the six points
blown up to obtain $X$:

\begin{cor}\label{typesvsneg}
Let $A$ and $A'$ be sets of six distinct points
of \pr2. Then $A$ and $A'$ have the same configuration type if and only if,
for some orderings $A=\{p_1,\ldots,p_6\}$ and $A'=\{p'_1,\ldots,p'_6\}$,
we have $f(\hbox{neg}(X))=\hbox{neg}(X')$, where $X$ is the surface
obtained by blowing up the points $p_i$, $X'$ is the surface obtained by blowing up
the points $p'_i$, $E_0,\ldots,E_6$ and $E'_0,\ldots,E'_6$ are the 
corresponding exceptional configurations
and $f:\hbox{Cl}(X)\to \hbox{Cl}(X')$ is the map defined
by $f(E_i)=E'_i$ for all $i$.
\end{cor}

\noindent {\em Proof}. 
If $A$ and $A'$ have the same configuration type,
then $f(\hbox{EFF}(X))=\hbox{EFF}(X')$,
hence $f(\hbox{NEG}(X))=\hbox{NEG}(X')$
(since \hbox{NEG} is the set of all $C$ in \hbox{EFF} such
that $C^2<0$ but $C$ is not the sum of two nontrivial elements of \hbox{EFF}),
so $f(\hbox{neg}(X))=\hbox{neg}(X')$. Conversely, by Remark \ref{minustwoRem},
$\hbox{neg}(X)$ determines $\hbox{NEG}(X)$, and, by Lemma \ref{NEGisfinite}
(and the proof of Lemma \ref{NEGisfinite}(c)), $\hbox{NEG}(X)$ determines
$h^0(X, G)$ for any class $G$. I.e., if  $f(\hbox{neg}(X))=\hbox{neg}(X')$,
then $h^0(X, G)=h^0(X', f(G))$ for every class $G$, hence $A$ and $A'$ have
the same configuration type. \prfend

The next remark shows explicitly how to determine Hilbert functions,
given $\hbox{NEG}(X)$ (or, equivalently by Remark \ref{minustwoRem}, given $\hbox{neg}(X)$).

\begin{rem}\label{HilbFuncRem}\rm
Given a fat points subscheme $Z=m_1p_1+\cdots+m_6p_6$ 
with support at 6 distinct points, for each $t$ consider
the class $F=F(Z,t)=tE_0-m_1E_1-\cdots-m_6E_6$. For each
$C\in \hbox{NEG}(X)$, check $F\cdot C$.
If $F\cdot C<0$, then $h_Z(t)=h^0(X, F)=h^0(X, F-C)$, so
we can replace $F$ by $F-C$ while preserving $h^0$.
Continue replacing the current $F$ by $F-C$ whenever the current $F$
meets some $C\in \hbox{NEG}(X)$ negatively. Eventually 
we obtain an $F$ such that either $F\cdot E_0<0$,
in which case $0=h^0(X, F)=h^0(X,F(Z,t))$, or
$F\cdot C\ge 0$ for all $C\in \hbox{NEG}(X)$, in which case
$F$ is nef and hence $h^0(X,F(Z,t))=h^0(X,F)$ is given
by Lemma \ref{NEGisfinite}(b). This procedure thus gives us a way to determine
the value $h_Z(t)$ of the Hilbert function $h_Z$
for every $t$. Note that determining $h_Z(t)$ involves
nothing more than integer arithmetic and addition and subtraction in
the rank 7 free abelian group $\hbox{Cl}(X)$. It
requires only that
we know $\hbox{NEG}(X)$ (or even just $\hbox{neg}(X)$) 
and the multiplicities $m_i$ of the
points of support of $Z$. We do not need to know the
points $p_i$ themselves. 

When $t\ge\alpha(Z)$, we also want to know
the multiplicity $m'_i=\hbox{mult}_{p_i}(C_{Z,t})$
and degree $d_{Z,t}$ of the curve $C_{Z,t}$ defined
by $\gamma(Z,t)$, whenever $\gamma(Z,t)$ has positive degree.
But $\gamma(Z,t)$ by Lemma \ref{NEGisfinite}
just defines the fixed component of the linear
system $I(Z)_t=H^0(X,F(Z,t))$, and hence
if $F$ is the nef divisor class obtained by successively
subtracting from $F(Z,t)$ classes in $\hbox{NEG}(X)$ as above,
then $F=F(Z_t^+,t-d_{Z,t})$ and 
$F-F(Z_t^+,t-d_{Z,t})=d_{Z,t}E_0-m'_1E_1-\cdots-m'_6E_6$,
so knowing $\hbox{NEG}(X)$ allows us to determine
$d_{Z,t}$ and the $m'_i$, and $Z_t^+$.
\end{rem}

Although Lemma \ref{NEGisfinite} gives us a criterion for a class being 
nef, our method of proof for Theorem \ref{distpntsThm} 
requires explicit generators for the nef cone; i.e., 
for the cone $\hbox{NEF}(X)$ of nef divisor classes on a given $X$, 
which by Lemma \ref{NEGisfinite} is just the cone 
of all $F$ such that $F\cdot C\ge 0$ for all $C\in\hbox{NEG}(X)$.
Actually, it will turn out that we will need explicit generators only when the anticanonical class,
$-K_X=3E_0-E_1-\cdots-E_6$, is nef. The problem of determining
generators of $\hbox{NEF}(X)$ is an example of the general
problem of finding generators for the dual of a nonnegative
subsemigroup whose generators are given
(in this case $\hbox{EFF}(X)$ is the subsemigroup, generated
by $\hbox{NEG}(X)$). This is not an easy computation in general, 
but in case $-K_X$ is nef the action of the
Weyl group, which we now recall, provides a significant simplification.

Let $r_0=E_0-E_1-E_2-E_3$ and for $1\le i\le 5$, let $r_i=E_i-E_{i+1}$.
(These are the so-called {\it simple\/} roots of the Lie-theoretic root system
of type ${\bf E}_6$.) Each homomorphism
$s_i:\hbox{Cl}(X)\to\hbox{Cl}(X)$ defined for any $x\in\hbox{Cl}(X)$
by the so-called {\it reflection\/} $s_i(x)=x+(x\cdot r_i)r_i$ through $r_i$
preserves the intersection product, and moreover $s_i(K_X)=K_X$
for all $i$. 
The subgroup of the orthogonal group of $\hbox{Cl}(X)$
generated by the $s_i$ is called the Weyl group, denoted $W_6$.
Since the reflection $s_i$ for $i>0$ is just the transposition
of $E_i$ and $E_{i+1}$, we see that $W_6$ contains the group
$S_6$ of permutations of $E_1,\ldots,E_6$. The element $s_0$ corresponds to
a quadratic transformation.

The group $W_6$ is a finite group of order 51,840. 
The $W_6$ orbit of $E_0$ is the following list
and those obtained from these by permuting the 
terms involving $E_i$ with $i>0$:\par
\tmprow{$E_0$}{$4E_0 -2E_1 -2E_2 -2E_3 -E_4 -E_5 -E_6$}
\tmprow{$2E_0 -E_1 -E_2 -E_3$}{$5E_0 -2E_1 -2E_2 -2E_3 -2E_4 -2E_5 -2E_6$}
\tmprow{$3E_0 -2E_1 -E_2 -E_3 -E_4 -E_5$}{}
\noindent Similarly, the $W_6$ orbit of $E_0-E_1$, up to permutations, is:\par
\tmprow{$E_0 -E_1$}{$3E_0 -2E_1 -E_2 -E_3 -E_4 -E_5 -E_6$}
\tmprow{$2E_0 -E_1 -E_2 -E_3 -E_4$}{}
\noindent And the $W_6$ orbit of $2E_0-E_1-E_2$, up to permutations, is:\par
\tmprow{$2E_0 -E_1 -E_2$}{$4E_0 -3E_1 -E_2 -E_3 -E_4 -E_5 -E_6$}
\tmprow{$3E_0 -2E_1 -E_2 -E_3 -E_4$}{$5E_0 -3E_1 -2E_2 -2E_3 -2E_4 -E_5 -E_6$}
\tmprow{$4E_0 -2E_1 -2E_2 -2E_3 -E_4 -E_5$}{$6E_0 -3E_1 -3E_2 -2E_3 -2E_4 -2E_5 -2E_6$}
\noindent The $W_6$ orbit of $3E_0 -E_1 -E_2 -E_3$, up to permutations, is:\par
\tmprow{$3E_0 -E_1 -E_2 -E_3$}{$6E_0 -4E_1 -2E_2 -2E_3 -2E_4 -E_5 -E_6$}
\tmprow{$4E_0 -2E_1 -2E_2 -E_3 -E_4$}{$7E_0 -4E_1 -3E_2 -3E_3 -2E_4 -2E_5 -E_6$}
\tmprow{$5E_0 -3E_1 -2E_2 -2E_3 -E_4 -E_5$}{$8E_0 -4E_1 -4E_2 -3E_3 -3E_4 -2E_5 -2E_6$}
\tmprow{$6E_0 -3E_1 -3E_2 -2E_3 -2E_4 -2E_5$}{$9E_0 -4E_1 -4E_2 -4E_3 -3E_4 -3E_5 -3E_6$}
\tmprow{$6E_0 -3E_1 -3E_2 -3E_3 -E_4 -E_5 -E_6$}{}
\noindent The $W_6$ orbit of $3E_0 -E_1 -E_2 -E_3-E_4$, up to permutations, is:\par
\tmprow{$3E_0 -E_1 -E_2 -E_3 -E_4$}{$5E_0 -3E_1 -2E_2 -2E_3 -E_4 -E_5 -E_6$}
\tmprow{$4E_0 -2E_1 -2E_2 -E_3 -E_4 -E_5$}{$6E_0 -3E_1 -3E_2 -2E_3 -2E_4 -2E_5 -E_6$}
\tmprow{$5E_0 -2E_1 -2E_2 -2E_3 -2E_4 -2E_5$}{$7E_0 -3E_1 -3E_2 -3E_3 -3E_4 -2E_5 -2E_6$}
\noindent The $W_6$ orbit of $3E_0 -E_1 -E_2 -E_3-E_4-E_5$, up to permutations, is:\par
\tmprow{$3E_0 -E_1 -E_2 -E_3 -E_4 -E_5$}{$5E_0 -2E_1 -2E_2 -2E_3 -2E_4 -2E_5 -E_6$}
\tmprow{$4E_0 -2E_1 -2E_2 -E_3 -E_4 -E_5 -E_6$}{}
\noindent Finally, the $W_6$ orbit of $-K_X=3E_0 -E_1 -E_2 -E_3-E_4-E_5-E_6$ is
just itself. The union of these orbits contains 1279 elements.
The next lemma says that the nef elements among
these 1279 generate the nef cone.

\begin{lem}\label{nefgenlem}
Let $X$ be a smooth projective
rational surface with a birational morphism to \pr2 such that
$\hbox{Cl}(X)$ has rank 7.
If $-K_X$ is nef, then 
the set $\Omega=\{F\in W_6G: F\cdot C\ge 0 \hbox{ for all }
C\in \C N\}$
generates $\hbox{NEF}(X)$ as a nonnegative subsemigroup
of $\hbox{Cl}(X)$, where 
$\C N$ is the set of classes
of reduced irreducible curves with $C^2=-2$ (the so-called {\it nodal roots})
and $G$ is the set consisting of
$E_0$,
$E_0 -E_1$,
$2E_0 -E_1 -E_2$,
$3E_0 -E_1 -E_2 -E_3$,
$3E_0 -E_1 -E_2 -E_3 -E_4$,
$3E_0 -E_1 -E_2 -E_3 -E_4 -E_5$, and
$3E_0 -E_1 -E_2 -E_3 -E_4 -E_5-E_6$.
\end{lem}

\noindent {\em Proof}. From the proof of Lemma \ref{NEGisfinite}, we know the complete
list of classes $C$ with $C^2=-2$ and $C\cdot K_X=0$
and it is not hard to check that they are contained in
(and, since $W_6$ preserves the intersection form, thus equal to) 
a single orbit of $W_6$; note, for example,
$s_3s_0(E_3-E_4)=E_0-E_1-E_2 -E_3$. 
This orbit is also known as the set of roots of the root system ${\bf E}_6$.
It is easy to verify that half of the roots are nonnegative integer linear 
combinations of the simple roots 
$r_0,\ldots, r_5$; the rest are the additive inverses of these.
The former are called {\it positive\/} roots; the latter are called 
{\it negative\/} roots. The class of a reduced irreducible curve $C$
with $C^2=-2$ is necessarily a positive root: it satisfies 
$C^2=-2$ and $C\cdot K_X=0$, so it is a root. Also, since 
$E_0$ is nef, we have $E_0\cdot C\ge0$. If $E_0\cdot C>0$, $C$ 
is clearly one of the positive roots. If $E_0\cdot C=0$, then
$C$ is a component of one of the exceptional curves $E_i$,
and thus of the form $E_i-E_j$ for some $0<i<j$, which is
a positive root. It is now not hard to check for any two
positive roots that $r\cdot r'\ge -2$, with $r\cdot r'=-2$ if and only if
$r=r'$.

Similarly, we also know the complete
list of classes $C$ with $C^2=-1$ and $C\cdot K_X=-1$, and we can
again check directly that they
form a single orbit $\C E$ of $W_6$; note, for example,
$s_0(E_1)=E_0 -E_2 -E_3$. Since $\C E$ is preserved under
the action of $W_6$, so is the nonnegative subsemigroup $\C E^*$
dual to $\C E$, consisting
of all classes $F$ such that $F\cdot C\ge 0$ for all $C\in\C E$.

By direct check, $G\cdot C\ge0$ for all $C\in\C E$,
so we have $G\subset{\C E}^*$, hence $W_6G\subset{\C E}^*$.
Since $\hbox{NEG}(X)=\C E\cup \C N$, it follows that $\Omega\subset\hbox{NEF}(X)$.
Now we must see that $\Omega$ generates $\hbox{NEF}(X)$.
Note that $\Omega$ is $W_6G\cap{\cal E}^*\cap{\cal N}^*$,
hence it is precisely the set of nef elements in
$W_6G$.

Since $W_6$ is finite, for each $F\in \C E^*$ there is some $w\in W_6$
such that $E_0\cdot wF$ is as small as possible. Let
$wF=a_0E_0-a_1E_1-\cdots-a_6E_6$. Since we can permute
the $a_i$ with $i>0$ by applying $s_j$ with $j>0$ and this
does not affect $E_0\cdot wF$, we may
assume that $a_1\ge a_2\ge \cdots\ge a_6$. Since $E_6\cdot wF\ge0$,
we have $a_6\ge0$. If $r_0\cdot wF<0$, then
we would have $E_0\cdot s_0(wF)<E_0\cdot wF$, so we 
also have $r_0\cdot wF\ge0$; i.e., $a_0\ge a_1+a_2+a_3$.

This means the $W_6$-orbit of every class $F\in  \C E^*$
intersects the subsemigroup $A$ of classes $H=b_0E_0-b_1E_1-\cdots-b_6E_6$
defined by the conditions $b_0\ge b_1+b_2+b_3$ and
$b_1\ge \cdots\ge b_6\ge 0$; i.e., by the conditions
$H\cdot r_0\ge0$, $\ldots$, $H\cdot r_5\ge0$.
It is not hard to check that
the set $G$ of classes $E_0$,
$E_0 -E_1$,
$2E_0 -E_1 -E_2$,
$3E_0 -E_1 -E_2 -E_3$,
$3E_0 -E_1 -E_2 -E_3 -E_4$,
$3E_0 -E_1 -E_2 -E_3 -E_4 -E_5$,
$3E_0 -E_1 -E_2 -E_3 -E_4 -E_5-E_6$
generates $A$ (which in fact turns out to be a fundamental
domain for the action of $W_6$ on $\C E^*$).
It is easy to check directly that, for every class $F$ in $A$,
$F\cdot C\ge0$ for every class $C$ with $C^2=-1$ 
and $C\cdot K_X=-1$. Also, $F\cdot r_i\ge 0$ holds for all $i$
since $F\in A$, hence $F\cdot C\ge0$ for the class $C$ of every reduced irreducible
curve with $C^2=-2$, since each such $C$ is a positive root.
Thus $A\subset \hbox{NEF}(X)$. 

Now let $F$ be any nef class. 
There is a sequence
$r_{i_1},\ldots,r_{i_l}$ of simple roots such that
$F_j\cdot C\ge0$ for all $C\in \C N_j$, each element of 
$\C N_j$ is a positive root, and  $F_l\in A$, where
$F=F_0$, $F_j=s_{i_j}(F_{j-1})$ for $1\le j\le l$,
$\C N_0=\C N$, and $\C N_j=s_{i_j}(\C N_{j-1})$ for $1\le j\le l$.
For each $j$, let $i_j$ be the largest $i$ such that
$F_{j-1}\cdot r_i<0$. If none exist, then
$l=j-1$ and $F_l\in A$, by definition of $A$.
Otherwise, let $F_j=s_{i_j}(F_{j-1})$. 
If $F=a_0E_0-a_1E_1-\cdots-a_6E_6$, what the sequence of 
operations does is to permute $a_1,\cdots,a_6$ so that
they are nondecreasing, and then to decrease $a_0$
whenever $s_0$ is applied. But the orbit $W_6F$ of $F$
is contained in $\C E^*$, hence every element
$H$ of the orbit has $H\cdot E_0=
H\cdot (E_0-E_1-E_2) + H\cdot E_1+H\cdot E_2\ge0$;
thus we cannot forever go on reducing the coefficient
of $E_0$, so eventually we arrive at a class $F_l$ for which
$F_l\cdot r_i\ge0$ for all $i$, and hence $F_l\in A$.
Now, $F_0\cdot C\ge0$ for all $C\in \C N_0$ since
$F=F_0$ is nef. Also, $wF\cdot wC=F\cdot C$ for all $w\in W_6$
since $W_6$ preserves the intersection form.
It follows that $F_j\cdot C\ge0$ for all $C\in \C N_j$ for all $j$.
Moreover, $r_{i_j}$ is never an element of
$\C N_{j-1}$, since $F_{j-1}\cdot r_{i_j}<0$. It is easy to check directly
that reflection by a simple root $r$ takes every positive root $r'\ne r$ to
another positive root. Thus each element of 
$\C N_j$ is a positive root for each $j$.

Since $F_l\in A$, $F_l$ is a nonnegative integer linear combination
of the classes in $G$. 
Moreover, the intersection of each of these classes with every element
of $\C N_l$ is nonnegative, since every element of $\C N_l$
is a positive root. Now let $w=s_{i_l}\cdots s_{i_1}$;
then $w^{-1}F_l=F$ and $w^{-1}H$ 
meets every element of $w^{-1}\C N_l=\C N$
nonnegatively for each $H\in G$. Thus each $w^{-1}H$ is nef, hence
$F$ is an integer linear combination of nef elements
in $W_6G$, as claimed.
\prfend

Given a nef divisor $F$, we still need a way of verifying that 
$\mu_F$ has maximal rank. Our main tools for doing so 
involve quantities 
$$q(F)=h^{0}({X},{F-E_1})\hbox{ and }l(F)=h^{0}({X},{F-(E_0-E_1})),$$ and
bounds on the dimension of the cokernel of
$\mu_F$, defined in terms of quantities
$q^*(F)=h^{1}({X},{F-E_1})$ and $l^*(F)=h^{1}({X},{F-(E_0-E_1)})$,
introduced in \cite{refIGC} and \cite{refFHH}. The 
following result is Lemma 2.2 of \cite{refFHH}. (There it is assumed
that $F\cdot E_1\ge F\cdot E_i$ for all $i>1$, but that is not needed
in the proof.)

\begin{lem}\label{IGClem}
Let $X$ be obtained by blowing up
distinct points $p_i\in\pr2$, and let
$F$ be the class of an effective divisor on $X$ with $h^1(X,F)=0$. Then 
$\hbox{dim ker $\mu_F$}\le q(F)+l(F)$
and $\hbox{dim cok $\mu_F$}\le q^*(F)+l^*(F)$.
\end{lem}

\begin{rem}\label{IGCrem}\rm
The quantities $q(F)$, $l(F)$, $q^*(F)$ and $l^*(F)$ are defined in terms of
$E_1$ and $E_0-E_1$, but in fact $E_j$, $j>0$, can
be used in place of $j=1$, since one can reindex the points.
\end{rem}

\begin{cor}\label{IGCcor}
Let $F$ and $G$ be nef divisors on
a surface $X$ obtained by blowing up
6 distinct points of \pr2.
If $q(F)>0$, $l(F)>0$ and $q^*(F)+l^*(F)=0$, 
then \hbox{dim cok $\mu_{F+G}=0$}.
\end{cor}

\noindent {\em Proof}. If more than three points are on a line,
then the six points are contained in a conic, and the result follows by
Theorem 3.1.2 of \cite{refFreeRes}. If at most
three points lie on a line, then, since there
are at most six points and they are distinct, $-K_X$ is nef.
So now we may assume $-K_X$ is nef.

That $q(F)>0$ implies $q(F+G)>0$
and $l(F)>0$ implies $l(F+G)>0$, are clear,
since a sum of effective divisors is effective.
By Lemma \ref{NEGisfinite}, $G+F$ is effective
and $h^1(X,G+F)=0$, so by Lemma \ref{IGClem} we have 
$\hbox{dim cok $\mu_{F+G}$}\le q^*(F+G)+l^*(F+G)$.
Thus it's enough to show $q^*(F+G)=0$ and $l^*(F+G)=0$.
By a direct check of the generators listed by Lemma \ref{nefgenlem}, $G$ 
is a sum of prime divisors of arithmetic genus at most 1. 
Hence it is enough by induction to show 
$q^*(F+G)=0$ and $l^*(F+G)=0$ when $G$ is the class of 
such a curve A. But this follows from
$0\to \C O_X(F-C)\to \C O_X(G+F-C)\to \C O_A(G+F-C)\to 0$,
taking $C$ to be $E_1$ (for $q^*$) or $E_0-E_1$ (for $l^*$),
since $h^1(X,F-C)=0$ by hypothesis, and $h^1(A,G+F-C)=0$.
(We have $A\cdot (G+F-C)\ge0$ since $G$ is nef, 
hence $h^1(A,G+F-C)=0$ if $A$ has genus 0, while
$G^2>0$ holds in each case that $A$ has genus 1. 
Thus $A\cdot (G+F-C)>0$ when the genus is 1,
hence again $h^1(A,G+F-C)=0$.)
\prfend

Given a nef divisor $F$, Corollary \ref{IGCcor} often applies,
in which case $\mu_{F+G}$ is surjective for all
nef $G$. However, not every nef class is an appropriate
sum of the form $F+G$. 
In the situations that we will need to deal with,
the set of those classes which are not of the appropriate form
turns out to be the union
of a finite set of exceptions (which we can handle by brute force)
with sets of strings of the form $F+iC$ (which 
we can handle by induction on $i$).

In order to set up the machinery to carry out the induction, 
define $\Gamma(X)$ to be the set of all nef classes 
which are not the sum of two nonzero nef classes. Then $\Gamma(X)$
generates $\hbox{NEF}(X)$ as a subsemigroup (i.e.,
every element of $\hbox{NEF}(X)$ is a nonnegative integer linear 
combination of elements of $\Gamma(X)$). 
For $i>0$, let $\Gamma_i(X)$ be the set of
all sums with exactly $i$ terms, where each term
is an element of $\Gamma(X)$. 
(So, for example, $\Gamma_1(X)=\Gamma(X)$.) Let
$S(X)$ be the set of all nef classes $F$ such that either 
$q(F)=0$, $l(F)=0$ or $l^*(F)+q^*(F)>0$.
Then let $S_i(X)=S(X)\cap\Gamma_i(X)$. By
Corollary \ref{IGCcor} we have $S_{i+1}(X)\subset S_i(X)+S_1(X)$.

Thus to show $\mu_F$ has maximal rank for 
every nef class $F$, it is enough by Lemma \ref{IGClem} to show that
$\mu_F$ has maximal rank for all $F\in S_i(X)$ for each $i$.
One checks directly that $\mu_F$ has maximal rank for all 
$F\in S_i(X)$ for small values of $i$. (It turns out
that it is never necessary to do this for $i>5$.)
For larger values of $i$, one applies Lemma \ref{stablem}
(the value of $k$ in this lemma never ends up needing to be bigger
than 2, although this is not obvious until after the fact) 
and Lemma \ref{stabindlem}. Also, it turns out that the inclusions
$S_{j+i}(X)\subset\{F+iC_F : F\in S_j(X)\}$ in 
Lemma \ref{stablem} can be chosen to be equalities, but that is more than we will need.

\begin{lem}\label{stablem}
Suppose for some $j$ there exists a $k$ and for each $F\in S_j(X)$
a $C_F\in S_1(X)$ such that
$S_{j+i}(X)\subset\{F+iC_F : F\in S_j(X)\}$ for $0\le i\le k$ and such that
whenever $C\in S_1(X)$ but $C\ne C_F$, then $F+kC\not\in S_{j+k}(X)$.
Then $S_{j+i}(X)\subset\{F+iC_F : F\in S_j(X)\}$
holds for all $i\ge0$.
\end{lem}

\noindent {\em Proof}. By Corollary \ref{IGCcor}, if $F+kC\not\in 
S_{j+k}(X)$, then $F+(k+1)C\not\in S_{j+k+1}(X)$.
Thus it is enough by induction to show 
$S_{j+k+1}(X)\subset\{F+(k+1)C_F : F\in S_j(X)\}$.
Say $G'\in S_{j+k+1}(X)$. Then $G'=G+C$, where
$G\in S_{j+k}(X)$ and $C\in S_{1}(X)$. By hypothesis,
$G=F'+kC_{F'}$ for some $F'\in S_j(X)$ and $C_{F'}\in S_1(X)$.
Since $G+C\in S_{j+k+1}(X)$, it follows by Corollary \ref{IGCcor} that
$F'+C\in S_{j+1}(X)$. Let $H=F'+C$; then $H=H'+C_{H'}$ for some
$H'\in S_j(X)$ and $C_{H'}\in S_1(X)$. 
Now, $H'+kC_{F'}\in S_{j+k}(X)$ (since $H+kC_{F'}=G+C\in S_{j+k+1}(X)$),
but for $D\in S_1(X)$ we have by hypothesis that $H'+kD\not\in S_{j+k}(X)$
unless $D=C_{H'}$. Thus $C_{F'}=C_{H'}$, so
$G+C=H'+(k+1)C_{H'}\in \{F+(k+1)C_F : F\in S_{j}(X)\}$,
so $S_{j+k+1}(X)\subset\{F+(k+1)C_F : F\in S_j(X)\}$.
\prfend

\begin{lem}\label{stabindlem}
Let $X$ be a blow up of \pr2 at
6 distinct points. Let $F$ be a nef divisor such that
$\mu_F$ is surjective, and let $C\subset X$ be the class of a smooth rational
curve such that $C^2\ge0$ and $(F+C)\cdot C\ge \hbox{max}(C\cdot E_1, C\cdot(E_0-E_1))$.
Then $\mu_{F+C}$ is surjective.
\end{lem}

\noindent {\em Proof}. Let $\Lambda$ denote $H^0(X,E_0)$, and apply the snake lemma to:

{
$$\matrix{
0 & \to & H^0(X,F)\otimes\Lambda & \to & H^0(X,F+C)\otimes\Lambda
& \to & H^0(C,\C O_X(F+C)\otimes\C O_C)\otimes\Lambda & \to & 0 \cr
  &     & \downarrow \mu_1 &     & \downarrow \mu_2 & & \downarrow \mu_3 & &  \cr
0 & \to & H^0(X,F+E_0) & \to & H^0(X,F+C+E_0) 
& \to & H^0(C,\C O_X(F+C+E_0)\otimes\C O_C) & \to & 0 \cr}$$
}

\noindent Since $\mu_F=\mu_1$ is onto, it is enough to show $\mu_3$ is
onto also, for which we apply 
$(F+C)\cdot C\ge \hbox{max}(C\cdot E_1, C\cdot(E_0-E_1))$,
using the criterion given in \cite{refFib}
(note also \cite{refFic}).
\prfend

We will be interested mostly in those $X$ such that
$2E_0-E_1-\cdots-E_6$ is not the class
of an effective divisor, since otherwise
(i.e., when the points $p_i$ lie
on a conic, possibly reducible or nonreduced)
$\mu_F$ is surjective whenever $F$ is nef
by Theorem 3.1.2 of \cite{refFreeRes},
which in turn depends on Lemma 2.5 of \cite{refFreeRes}.
However, some details were left out of the published
proof of this lemma, so we present it here in full.
The extra details are indicated by indentation.

\begin{lem}\label{keylemma}
Let $X$ be a smooth projective 
rational surface, and let \C N be the class of 
an effective divisor $N$ on $X$
such that $h^0(X, \C N+K_X)=0$. If \C F and \C G 
are the restrictions to $N$
of divisor classes $\C F'$ and $\C G'$ on $X$ 
which meet each component of $N$ nonnegatively,
then $\C S(\C F,\C G)=0$, where $\C S(\C F,\C G)$ denotes
the cokernel of the natural map
$H^0(N,{\cal F})\otimes H^0(N,{\cal G})\to H^0(N,{\cal F}+{\cal G})$.
\end{lem}

\noindent {\em Proof}. To prove the lemma, induct on the sum $n$ of 
the multiplicities of the components of $N$.
By Lemma II.9 of \cite{refARS}, $h^1(N, \C O_N)=0$ and every component
of $N$ is a smooth rational curve. Thus the case $n=1$
is trivial (since then $N=\pr 1$, and
the space of polynomials of degree $f$ in two variables 
tensor the space of polynomials of degree $g$ in two variables 
maps onto the space of polynomials of degree $f+g$). So say $n>1$.

As in the proof
of Theorem 1.7 of \cite{refA}, $N$ has a component $C$ such that 
$(N-C)\cdot C\le 1$. Let $L$ be the effective divisor $N-C$
and let \C L be its class.
Thus we have an exact sequence 
$0\to \C O_C\otimes(-\C L)\to \C O_N\to \C O_L\to 0$.

\begin{quote}
\ \hbox to20pt{\hfil}To see this, apply the snake lemma to 
$$\matrix{
0 & \to & \C O_X(-N) & \to & \C O_X & \to & \C O_N & \to & 0 \cr
  &     & \downarrow &     & \downarrow & & \downarrow & & \cr
0 & \to & \C O_X(-L) & \to & \C O_X & \to & \C O_L & \to & 0 \cr}
$$
to see that the kernel of $\C O_N \to \C O_L$ is just
the cokernel of $\C O_X(-N) \to \C O_X(-L)$, which is just
$\C O_C\otimes \C O_X(-L)$, which we may write as $\C O_C(-L)$.
\end{quote}

Now, $-L\cdot C\ge -1$, and both 
$\C F'$ and $\C G'$ meet $C$ nonnegatively. We may assume
$\C F'\cdot C\ge\C G'\cdot C$, otherwise reverse the roles of 
$\C F'$ and $\C G'$. Since
$C=\pr 1$, we see that $h^1(C, \C O_C\otimes(\C F'-\C L))$,
$h^1(C, \C O_C\otimes(\C G'-\C L))$ and 
$h^1(C, \C O_C\otimes(\C F'+\C G'-\C L))$ all vanish.
An argument similar to that used to prove Proposition II.3(a, b) 
of \cite{refFreeRes}
now shows that we have an exact sequence
$\C S(\C O_C\otimes(\C F'-\C L),\C O_C\otimes\C G')\to
\C S(\C F,\C G)\to\C S(\C O_L\otimes\C F,\C O_L\otimes\C G)\to0$.

\begin{quote}
\ \hbox to20pt{\hfil}What is actually clear here is that we have
$\C S(\C O_C\otimes(\C F'-\C L),\C G)\to
\C S(\C F,\C G)\to\C S(\C O_L\otimes\C F,\C G)\to0$.
Since $h^1(C, \C O_C\otimes(\C G'-\C L))=0$, we know
$\C G\to \C O_L\otimes\C G$ is surjective on global sections,
and hence that $\C S(\C O_L\otimes\C F,\C G)$ is the same as
$\C S(\C O_L\otimes\C F,\C O_L\otimes\C G)$.
What needs additional justification here is that
$\C O_N\otimes\C G'\to \C O_C\otimes\C G'$ is surjective on global sections,
so that we can conclude that 
$\C S(\C O_C\otimes(\C F'-\C L),\C G)$ is the same as
$\C S(\C O_C\otimes(\C F'-\C L),\C O_C\otimes\C G')$.
             
\ \hbox to20pt{\hfil}Now, $\C N+K_X$ is not the class of an effective divisor,
and the same will remain true if we replace $N$ by any subscheme 
of $N$ obtained by subtracting off
irreducible components of $N$.
Thus any such resulting subscheme $M$ of $N$ has the property, like $N$ itself,
that there is a component $D$ of $M$ such that 
$(M-D)\cdot D\le 1$. If $M$ is just $N$ with the reduced induced
scheme structure, then by induction on the number of components
of $M$ it follows (using Lemma II.9 of \cite{refARS}) that any two components of $N$ are smooth
rational curves that are either disjoint or meet transversely
at a single point, and no sequence
$B_1$, $\ldots$, $B_i$ of distinct components exists
such that $B_i\cdot B_{1}>0$ and $B_j\cdot B_{j+1}>0$ for $1\le j<i$
(in particular, no three components meet at a single point,
and the components of $M$ form a disjoint union of trees).

\ \hbox to20pt{\hfil}First assume that $N$ is reduced; i.e. that $N=N_{red}$.
Then $C$ is not a component of $N-C$.
Choose a section $\sigma_C$ of $\C O_C\otimes\C G'$,
and for each of the other components $B$ of $N$,
choose a section $\sigma_B$ of $\C O_B\otimes\C G'$
such that $\sigma_B$ does not vanish at any of the points where
$B$ meets another component of $N$. (This is possible since
$B$ is smooth and rational, so $\C O_B\otimes\C G'$
is $\C O_{\pr1}(d)$ for some $d\ge 0$, so a section can always
be chosen which does not vanish at any of a given
finite set of points of $B$.) Since $N$ has no cycles
and the components meet transversely,
it is clear that starting from $\sigma_C$
one can patch together appropriate scalar multiples
of the sections $\sigma_B$
to get a section $\sigma$ of \C G which restricts
to $\sigma_C$. Thus $\C O_N\otimes\C G'\to \C O_C\otimes\C G'$
is surjective on global sections.

\ \hbox to20pt{\hfil}Now assume that $N$ is not reduced. Let $M$ be the union
of the components of $N$ which have multiplicity greater than 1 
(taken with the same
multiplicities as they have in $N$) together with those multiplicity 1 components
of $N$ that meet one of these. No multiplicity 1 component
$B$ of $M$ satisfies $B\cdot (M-B)\le 1$, so there must be a
component $B$ of multiplicity more than 1 that does, 
and hence we also have
$B\cdot (N-B)\le 1$ for some component $B$ of $N$ of 
multiplicity more than 1. Now from this and
$0\to \C O_B(-N+B)\otimes 
\C G\to \C O_N\otimes\C G'\to \C O_J\otimes\C G'\to0$,
where $J=N-B$, we see $h^1(B, \C O_B(-N+B)\otimes \C G) = 0$, so
$\C O_N\otimes\C G'\to \C O_J\otimes\C G'$ is surjective on global sections.
But $J$ still has $C$ as a component, 
because either $C$ has multiplicity 1 in $N$ (and hence
$C\ne B$), or $C$ has multiplicity more than 1 in $N$ (and so even if $B=C$, $C$ remains
a component of $N-B=J$). By induction on the number of components,
we conclude that $\C O_N\otimes\C G'\to \C O_{N_{red}}\otimes\C G'$
is surjective on global sections. But $C$ is still a component of
$N_{red}$, and $\C O_{N_{red}}\otimes\C G'\to \C O_{C}\otimes\C G'$
is surjective on global sections from above, hence so is 
$\C O_{N}\otimes\C G'\to \C O_{C}\otimes\C G'$.
\end{quote}

Since $\C S(\C O_L\otimes\C F,\C O_L\otimes\C G)=0$ by induction,
it suffices to show
$\C S(\C O_C\otimes(\C F'-\C L),\C O_C\otimes\C G')=0$.
If $C\cdot(\C F'-\C L)\ge 0$, then the latter is 0 (as in the previous paragraph).
Otherwise, we must have $0=\C F'\cdot C=\C G'\cdot C$ and
$C\cdot L=1$, so
$\C O_C(-1)=\C O_C\otimes(\C F'-\C L)$ and
$\C O_C=\C O_C\otimes\C G'$, which means 
$h^0(\C O_C,\C O_C\otimes(\C F'+\C G'-\C L))=0$ 
and hence again $\C S(\C O_C\otimes(\C F'-\C L),\C O_C\otimes\C G')=0$. \prfend

\section{The Main Results}\label{distpnts}

In this section we first determine, up to permuting
$E_1,\ldots,E_6$, which subsets of $\C L\cup\C Q$          
occur as subsets of the form $\hbox{neg}(X)$, which by Corollary \ref{typesvsneg}
is equivalent to determining the configuration types for six distinct points of \pr2.
What we find is that the types are precisely those shown in Figure 1, where the classes of the 
proper transforms of the curves shown in a diagram of Figure 1
give the elements of $\hbox{neg}(X)$ for the corresponding configuration type.
We then prove our main result, Theorem \ref{distpntsThm},
and finish by explicitly answering, in the case of 6 points, 
the questions raised in \cite{refGMS}.

To begin, note that the elements $C$ of $\hbox{neg}(X)$ satisfy
the following three conditions: (i) $C\in \C L\cup\C Q$; (ii) 
$C^2<-1$; and (iii) $C\cdot D\ge 0$ whenever $C,D\in \hbox{neg}(X)$ with $C\ne D$.

First, if $2E_0-(E_1+\cdots+E_6)\in\hbox{neg}(X)$,
then $\{2E_0-(E_1+\cdots+E_6)\}=\hbox{neg}(X)$.
(For if $C\in \hbox{neg}(X)$ but $C\ne 2E_0-(E_1+\cdots+E_6)$,
then $C\cdot (2E_0-(E_1+\cdots+E_6))\ge 0$ by (iii). But by direct check,
every element $C\in \C L\cup\C Q$ with $C^2<-1$ has
$C\cdot (2E_0-(E_1+\cdots+E_6))<0$.) 
The case that $\{2E_0-(E_1+\cdots+E_6)\}=\hbox{neg}(X)$
corresponds to configuration type 11 in Figure 1. It is clear that this possibility 
actually occurs, since blowing up any six points on a
smooth conic results in $2E_0-(E_1+\cdots+E_6)\in\hbox{neg}(X)$,
and hence, as we just saw, $\{2E_0-(E_1+\cdots+E_6)\}=\hbox{neg}(X)$.

We now classify sets $M$ satisfying the conditions:
(i)  $M\subset \C L$ ; (ii) if $C\in M$, then
$C^2<-1$; and (iii) $C\cdot D\ge 0$ whenever $C, D\in M$, $C\ne D$.
For each such $M$, we also will show that there is an $X$ with
$M=\hbox{neg}(X)$.

In fact, such a subset $M$ is just a matroid of rank 3 or less on a six point set,
or, in the terminology of \cite{refBCH}, it is a plane 6 point combinatorial geometry.
It is not hard to work them all out, but \cite{refBCH} gives a complete list,
saving us the trouble of doing so. The result corresponds precisely with
what we show as configuration types 1 through 10 in Figure 1.
So now we merely need to see that they all arise.

To show configuration type 1 occurs, we just need to show that
one can pick 6 points such that no line passes through
any 3 and no conic passes through all 6. Thus we can
pick any two distinct points to be $p_1$ and $p_2$.
Then $p_3$ can be any point not on the line through
$p_1$ and $p_2$; $p_4$ can be any point not on any line through
two of the first three points, and $p_5$ can be any point not on any line through
two of the first four points. Finally, $p_6$ can be any point not on any line through
two of the first five points nor on the conic through the first five points
(of which there is only one). At each step we are allowed to choose any point
avoiding a proper closed subset of \pr2. There is no obstruction to doing
this, so configuration type 1 occurs.

For configuration type 2, we proceed as before, but the last point
must be on exactly one of the lines through two of the
previously chosen points. For example, we choose $p_6$
to be on the line $L$ through $p_1$ and $p_2$, but not on
any other line through two of the previously chosen points.
Thus the condition on our choice of $p_6$ is that we avoid finitely
many points of $L$, which clearly we may do.

By similar reasoning, it is easy to check that each of the configurations
1 through 9 occur. With configuration 10, the same reasoning works 
to choose points $p_1$ through $p_5$, but the choice of $p_6$ is forced,
since $p_1,\ldots,p_5$ uniquely determine $p_6$. Since we have no
freedom in our choice of $p_6$, our previous argument is
invalid at the last step. Instead, we take
our six points to be the points of intersection of four general lines.
Clearly, no four of the points can be collinear. 
So now we must check that the four lines are the only lines
through any three of the points. Suppose there were a
fifth line going through three of the points. Given 
any three of the six points of intersection of four general lines,
it is easy to check that one of the four lines passes through two of the three points.
So there can be no fifth line through any three of the points.
Thus configuration type 10 also occurs.

(The foregoing justifications that the configurations actually occur
may at first sight seem unnecessary. To show that they are not, we mention 
a similar example involving seven points. Suppose we want a configuration
of six lines through seven points, arranged such that each line
passes through exactly 3 points. Intuitively, we get such a configuration
by taking three of the lines to be sides of an equilateral triangle, and the other three
to be the angle bisectors. The seven points are the points where any two
of the lines meet. This configuration occurs if and only if the ground field
does not have characteristic 2. When the  characteristic is 2,
an additional line through the midpoints on the sides of the triangle
is forced.)

We now prove our main result:

\begin{thm}\label{distpntsThm}
Let $X$ be obtained by blowing up 
6 distinct points of \pr2. Let $E_0$, $E_1$, $\ldots$, $E_6$ be the 
corresponding exceptional configuration.
Let $F$ be a nef divisor on $X$. Then $\mu_F$ has maximal rank.
\end{thm}

\noindent {\em Proof}. We first consider the two extremes.
If no line contains three of the points
and no conic contains all 6, 
then the result follows by \cite{refFi}.
This is the case in which the points are general.
If all 6 points are on a conic, the result follows by
Theorem 3.1.2 of \cite{refFreeRes}. Note also that
if 4 or more of the points are on a line, then
all 6 are on a conic.

So now we are reduced to considering the case that
some line contains three points, but no line
contains 4 or more of the points and no conic contains all 6. 
Thus $\hbox{neg}(X)$ consists only of classes 
of the form $L-E_{i_1}-E_{i_2}-E_{i_3}$, hence 
$\hbox{neg}(X) =\C N$.
If there is more than one line
that contains three of the points, then any two such lines
must share a point (otherwise all 6 points would lie on the
two lines, which is a conic). It follows that
 the set $\C N$ of nodal roots must, up to
indexation of the points $p_i$, be one of the following:
\begin{description}
\item[(i)] $\{r_0\}$ --- i.e., the first three points are on a line
and no other set of three points is on a line (this 
case corresponds to configuration type 2);
\item[(ii)] $\{r_0, E_0-E_1-E_4-E_5\}$ --- i.e., two of the points
are on one line, two on another, a fifth point
occurs where the two lines meet, and the sixth point
is not on any line through any two of the other points
(this case corresponds to type 8);
\item[(iii)] $\{r_0, E_0-E_1-E_4-E_5, E_0-E_3-E_5-E_6\}$ --- i.e., 
three lines form a triangle, with three of the points at the vertices,
with an additional point on each line, but these last three points
are not collinear (this case corresponds to type 9);
\item[(iv)] $\{r_0, E_0-E_1-E_4-E_5, E_0-E_3-E_5-E_6, E_0-E_2-E_4-E_6\}$ --- i.e., 
the 6 points are the points of intersection of four lines,
no three of which meet at a single point (this 
case corresponds to type 10).
\end{description}

We now treat case (iv) in detail. The other cases 
(and the case that $\C N$ is empty, which thereby recovers the result
of \cite{refFi}) are similar. Using Remark \ref{minustwoRem}, from $\C N=\hbox{neg}(X)$
we determine that $\hbox{NEG}(X)$ consists of the following
classes (where we list only the coefficients,
so, for example, \hbox{\tt 1  0 -1 0 -1 0 -1}
denotes $E_0-E_2-E_4-E_6$):

{\footnotesize
\parindent=30pt
\begin{verbatim}
                0  1  0  0  0  0  0     1  0  0 -1 -1  0  0     1 -1 -1 -1  0  0  0
                0  0  1  0  0  0  0     1  0 -1  0  0 -1  0     1 -1  0  0 -1 -1  0
                0  0  0  1  0  0  0     1 -1  0  0  0  0 -1     1  0  0 -1  0 -1 -1
                0  0  0  0  1  0  0                             1  0 -1  0 -1  0 -1
                0  0  0  0  0  1  0  
                0  0  0  0  0  0  1
\end{verbatim}}

Next we need to determine generators for $\hbox{NEF}(X)$.
By Lemma \ref{nefgenlem}, the set of all $F\in W_6G$
such that $F\cdot C\ge0$ for all 
$C\in \C N=\{r_0, E_0-E_1-E_4-E_5, E_0-E_3-E_5-E_6, E_0-E_2-E_4-E_6\}$
generates $\hbox{NEF}(X)$, where $W_6G$ is the set
of 1279 elements of the $W_6$ orbits
of the elements of $G$ from Lemma \ref{nefgenlem}. A tedious but easily
coded check results in 212 generators.
Many of these 212 are a sum of two other classes among
the 212. Removing all classes which occur as such sums, 
we are left with 39, which therefore generate.
Here is a list of these 39:

{\footnotesize
\parindent=30pt
\begin{verbatim}
                1  0  0  0  0  0  0     2 -1  0 -1  0 -1  0     3  0  0 -1 -2 -1 -1
                2  0 -1 -1 -1  0  0     2 -1  0  0 -1  0 -1     3 -1  0 -2 -1 -1  0
                2  0  0 -1 -1 -1  0     2 -1 -1  0  0 -1  0     3  0 -1  0 -1 -2 -1
                2  0  0  0 -1 -1 -1     2 -1  0 -1  0  0 -1     3 -1 -1 -1  0  0 -2
                2 -1  0 -1 -1  0  0     1 -1  0  0  0  0  0     3 -1 -1 -1 -2  0  0
                2  0 -1  0 -1 -1  0     1  0 -1  0  0  0  0     3 -1  0  0 -1 -1 -2
                2  0  0 -1 -1  0 -1     1  0  0 -1  0  0  0     3  0 -1 -2 -1  0 -1
                2  0 -1  0  0 -1 -1     1  0  0  0 -1  0  0     3 -1 -2  0 -1 -1  0
                2  0 -1 -1  0  0 -1     1  0  0  0  0 -1  0     3  0 -2 -1  0 -1 -1
                2 -1  0  0  0 -1 -1     1  0  0  0  0  0 -1     3 -1 -1 -1  0 -2  0
                2 -1 -1  0  0  0 -1     2  0 -1 -1 -1 -1  0     3 -2  0 -1  0 -1 -1
                2 -1 -1  0 -1  0  0     2 -1 -1  0  0 -1 -1     3 -2 -1  0 -1  0 -1
                2  0 -1 -1  0 -1  0     2 -1  0 -1 -1  0 -1     3 -1 -1 -1 -1 -1 -1
\end{verbatim}}

For each of these classes $F$, we find (after reindexing if need be,
as discussed in Remark \ref{IGCrem}, but using the same indexing for $q$, $l$,
$q^*$, and $l^*$) that $q^*(F)=0=l^*(F)$,
hence $\mu_F$ is surjective by Lemma \ref{IGClem} and
Remark \ref{IGCrem}. For example, to see how to 
compute these quantities, consider
$q^*(F)$ for $F=3E_0 -E_1 -2E_3 -E_4 -E_5$
from the list above. Then, applying Remark \ref{IGCrem},
we reindex so that $q(F)=h^0(X,F-E_3)$
and $l(F)=h^0(X,F-(E_0-E_3))$, etc.
Since $r_0\in \hbox{NEG}(X)$ and 
$r_0\cdot (F-E_3)<0$, we see 
$h^0(X,F-E_3)=h^0(X,F-E_3-r_0)$.
But now $E_2\cdot (F-E_3-r_0)<0$, so
now $h^0(X,F-E_3)=h^0(X,F-E_3-r_0-E_2)$.
Continuing in this way we eventually find
that  $h^0(X,F)=\cdots=h^0(X,0)=1$,
hence $q(F)=1$. Riemann-Roch now states that $q(F)-q^*(F)=
((F-E_3)^2-(F-E_3)\cdot K_X)/2+1=1$,
so $q^*(F)=0$.

Of the 39, all but the following 9 have both $q$ and $l$ 
positive, and thus $S_1(X)$ is just the set of these 9:

{\footnotesize
\parindent=30pt
\begin{verbatim}
                1 -1  0  0  0  0  0     1  0  0  0 -1  0  0     2  0 -1 -1 -1 -1  0
                1  0 -1  0  0  0  0     1  0  0  0  0 -1  0     2 -1 -1  0  0 -1 -1
                1  0  0 -1  0  0  0     1  0  0  0  0  0 -1     2 -1  0 -1 -1  0 -1
\end{verbatim}}

In each of these cases $q=0$. By Corollary \ref{IGCcor}, $\mu_F$
is surjective for all nef $F$ except possibly those in the
subsemigroup generated by these last 9. A direct check shows that 
the conditions of Lemma \ref{stablem} apply here with $k=2$ and $C_F=F$,
so $S_{i}(X)\subset\{iF: F\in S_1(X)\}$ for all $i$. Surjectivity for
$\mu_{iF}$ for each $F$ and $i$ follows by
direct check that $q^*(iF)=0=l^*(iF)$ when $i=1$ and 2, and
then for all $i>0$ by applying Lemma \ref{stabindlem}.

Cases (i), (ii) and (iii) are handled the same way, thereby proving 
Theorem \ref{distpntsThm}. For case (i), $S_1(X)$ has 55 elements,
$S_2(X)$ has 90 elements, and $S_i(X)$ has 93 elements
for $i>2$. Lemma \ref{stablem} applies for $j=3$ with $k=2$, although
this time it is not always true that $F$ is a multiple of $C_F$.
For example, $F=7E_0 -2E_1 -\cdots-2E_5 -5E_6\in S_3(X)$,
but $C_F=3E_0 -1E_1 -\cdots -1E_5 -2E_6$.
For case (ii), $S_1(X)$ has 37 elements, $S_i(X)$ has 34 elements
for $i>1$ and Lemma \ref{stablem} applies for $j=2$ with $k=2$.
For case (iii), $S_1(X)$ has 22 elements, $S_i(X)$ has 12 elements
for $i>1$ and Lemma \ref{stablem} applies for $j=2$ with $k=2$.
(For the case that $\C N$ is empty, $S_1(X)$ has 159 elements, $S_2(X)$ has 301 elements,
and $S_i(X)$ has 316 elements for $i>2$. Lemma \ref{stablem} applies for $j=3$ with $k=2$.
Lemma \ref{stabindlem} then gives the result except for multiples of
$F=5E_0-2E_1-\cdots-2E_6$, since $\mu_F$ is injective,
and $l^*(mF)=1$ for $m\ge0$. Thus one must show ad hoc that
$\mu_{2F}$ is surjective (see \cite{refFi}); then Lemma \ref{stabindlem}
applies to show that $\mu_{mF}$ is surjective for all $m>2$.)
\prfend

\begin{example}\label{distptsex}\rm
We work out an example to show how to determine the
Hilbert function and graded Betti numbers of the ideal of a fat point subscheme.
Assume the points are arranged as in case (iv); that is, configuration type 10.
Assume the points are indexed so that a line passes through
points 1, 2 and 3, and through 1, 4 and 5, and 2, 4 and 6 and
3, 5 and 6.
Let $Z=2p_1+2p_2+6p_3+2p_4+2p_5+2p_6$.
The associated divisor class for degree $i$ is
$F(Z,i)=iE_0-(2E_1+2E_2+6E_3+2E_4+2E_5+2E_6)$.
Computing $h_Z(i)=h^0(X, F(Z,i))$ as in Remark \ref{HilbFuncRem},
we find $h_Z(5)=0$, $h_Z(6)=1$, 
$h_Z(7)=4$, $h_Z(8)=11$, $h_Z(9)=19$ and $h_Z(10)=30$,
so $\alpha(Z)=6$.
Also,  $h^1(X, F(Z,8))>0$ but $h^1(X, F(Z,9))=0$,
hence the regularity $\sigma(Z)$ is 10. 

Thus $t_i=0$ for $i<\alpha(Z)=6$
and for $i>\sigma(Z)=10$, and since
$h_Z(6)=1$, we see $t_6=1$ and
that $\mu_{F(Z,6)}$ is injective so
$t_7=h_Z(7)-3h_Z(6)=1$. 
To find $t_8$, note that:
$F(Z,7)\cdot C_1<0$, where $C_1=E_0-E_3-E_4$;
$(F(Z,7)- C_1)\cdot C_2 < 0$ for $C_2=r_0$;
$(F(Z,7)- C_1-C_2)\cdot C_3 < 0$ for $C_3=E_0-E_3-E_4-E_5$;
$(F(Z,7)- C_1-C_2-C_3)\cdot C_2 < 0$;
$(F(Z,7)- C_1-2C_2-C_3)\cdot C_3 < 0$;
and $F(Z,7)- C_1-2C_2-2C_3$ is nef.
Thus the divisor class of fixed components of $F(Z,7)$ 
is $C_1+2C_2+2C_3=5E_0-2E_1-2E_2-5E_3-E_4-2E_5-2E_6$,
so $Z_7^-=2p_1+2p_2+5p_3+p_4+2p_5+2p_6$, $d_{Z,7}=5$,
and $Z_7^+=p_3+p_4$. Now we have
$t_8=\hbox{dim cok}(\mu_{F(Z,7)- C_1-2C_2-2C_3})+
(h^0(X,E_0+F(Z,7))-h^0(X,E_0+F(Z,7)- C_1-2C_2-2C_3))$.
But $F(Z,7)- C_1-2C_2-2C_3$ is nef, its $\mu$ is onto by Theorem \ref{distpntsThm},
and $h^0(X,E_0+F(Z,7))-h^0(X,E_0+F(Z,7)- C_1-2C_2-2C_3)=
h_Z(8)-h^0(X,E_0+F(Z,7)- C_1-2C_2-2C_3)=11-8=3$.
Similarly, $t_9=0$ and $t_{10}=2$.
From the triple difference $\Delta^3 h_Z$, we find
$s_i=0$ except for $s_8=1$, $s_9=3$ and $s_{11}=2$.
Thus the minimal free resolution of $I_Z$ is
$0\to F_1\to F_0\to I_Z\to 0$
where $F_0=R[-6]\oplus R[-7]\oplus R[-8]^3\oplus R[-10]^2$
and $F_1=R[-8]\oplus R[-9]^3\oplus R[-11]^2$.
\end{example}

It is easy to implement the procedure demonstrated
in Example \ref{distptsex} as, for example, an awk script. 
We did so; the resulting script can be run over the web by visiting

\noindent {\tt http://www.math.unl.edu/$\scriptstyle\sim$bharbour/6ptres/6reswebsite.html} .

\noindent We used it to determine the Hilbert functions and 
graded Betti numbers for the ideals defining 
$Z=p_1+\cdots+p_6$ and for $2Z=2p_1+\cdots+2p_6$
for each of the 11 configuration types, thereby answering in the case
of six points the questions raised in \cite{refGMS}. We could just as easily
run $mZ$ for any $m$ or for any multiplicities
$m_1p_1+\cdots+m_6p_6$, if we wished to answer
the questions raised by \cite{refGMS} not only for double points
but for points of any given multiplicities. 
Note that for configuration types 5, 7 and 11, $Z$ is a complete intersection,
and thus the Hilbert function and graded Betti numbers
for $mZ$ are already known for all $m$ (see, for example,
\cite{refBGVTone} and \cite{refBGVTtwo}). Also, the Hilbert function and graded Betti numbers
for $m_1p_1+\cdots+m_6p_6$ for any $m_i$ are known by
\cite{refFi} for configuration type 1, and by \cite{refFreeRes} for
configurations 3, 4, 5, 6, 7 and 11 (since the points are contained in a conic).
Results for configuration types 2, 8, 9 and 10 are new. 

For ease of comparison with
results of \cite{refGMS}, we give the Hilbert functions
$h_{R/I(Z)}$ of $R/I(Z)$, rather than for $I(Z)$. The Hilbert function
of $R/I(Z)$ in degree 0 is always 1, and then it
increases until it achieves the value $\hbox{deg}(Z)$,
at which point it becomes constant. In each case we show
the value $h_{R/I(Z)}(t)$ of the Hilbert function in each degree $t\ge0$
until it becomes constant.

Here are the results. There are four different Hilbert functions
for $Z$, and all together there are six different Hilbert functions
for $2Z$, two whose support has one Hilbert function,
two whose support has another, and one each for the remaining
two cases. Note that for each Hilbert function for $Z$,
there is among the Hilbert functions for $2Z$
both a maximum and minimum Hilbert function.
\vskip\baselineskip

{\footnotesize
\mytmprow{Scheme}{Type(s)}{$h_{R/I(Z)}$}{}
\vskip-5pt
\hbox to\hsize{\hrulefill}
\mytmprow{$Z$}{1, 2, 8, 9, 10}{1, 3, 6, 6}{}
\mytmprow{$2Z$}{1, 2, 8, 9}{1, 3, 6, 10, 15, 18, 18}{}
\mytmprow{$2Z$}{10}{1, 3, 6, 10, 14, 18, 18}{}
\vskip-5pt
\hbox to\hsize{\hrulefill}
\mytmprow{$Z$}{3, 6, 7, 11}{1, 3, 5, 6, 6}{}
\mytmprow{$2Z$}{3, 6}{1, 3, 6, 10, 14, 16, 17, 18, 18}{}
\mytmprow{$2Z$}{7, 11}{1, 3, 6, 10, 14, 17, 18, 18}{}
\vskip-5pt
\hbox to\hsize{\hrulefill}
\mytmprow{$Z$}{4}{1, 3, 4, 5, 6, 6}{}
\mytmprow{$2Z$}{4}{1, 3, 6, 10, 12, 14, 15, 16, 17, 18, 18}{}
\vskip-5pt
\hbox to\hsize{\hrulefill}
\mytmprow{$Z$}{5}{1, 2, 3, 4, 5, 6, 6}{}
\mytmprow{$2Z$}{5}{1, 3, 5, 7, 9, 11, 13, 14, 15, 16, 17, 18, 18}{}
\vskip-5pt
\hbox to\hsize{\hrulefill}

\myrow{Scheme}{Type(s)}{$F_1$}{$F_0$}
\vskip-5pt
\hbox to\hsize{\hrulefill}
\myrow{$Z$}{1, 2, 8, 9, 10}{$R[-4]^3$}{$R[-3]^4$}
\myrow{$2Z$}{1, 2}{$R[-7]^3$}{$R[-6]\oplus R[-5]^3$}
\myrow{$2Z$}{8}{$R[-7]^3\oplus R[-6]$}{$R[-6]^2\oplus R[-5]^3$}
\myrow{$2Z$}{9}{$R[-7]^3\oplus R[-6]^2$}{$R[-6]^3\oplus R[-5]^3$}
\myrow{$2Z$}{10}{$R[-7]^4$}{$R[-6]^4\oplus R[-4]$}
\vskip-5pt
\hbox to\hsize{\hrulefill}
\myrow{$Z$}{3, 6}{$R[-5]\oplus R[-4]$}{$R[-4]\oplus R[-3]\oplus R[-2]$}
\myrow{$2Z$}{3}{$R[-9]\oplus R[-7]\oplus R[-6]$}{$R[-8]\oplus R[-5]^2\oplus R[-4]$}
\myrow{$2Z$}{6}{$R[-9]\oplus R[-7]\oplus R[-6]^2$}{$R[-8]\oplus R[-6]\oplus R[-5]^2\oplus R[-4]$}
\vskip-5pt
\hbox to\hsize{\hrulefill}
\myrow{$Z$}{7, 11}{$R[-5]$}{$R[-3]\oplus R[-2]$}
\myrow{$2Z$}{7, 11}{$R[-8]\oplus R[-7]$}{$R[-6]\oplus R[-5]\oplus R[-4]$}
\vskip-5pt
\hbox to\hsize{\hrulefill}
\myrow{$Z$}{4}{$R[-6]\oplus R[-3]$}{$R[-5]\oplus R[-2]^2$}
\myrow{$2Z$}{4}{$R[-11]\oplus R[-7]\oplus R[-5]^2$}{$R[-10]\oplus R[-6]\oplus R[-4]^3$}
\vskip-5pt
\hbox to\hsize{\hrulefill}
\myrow{$Z$}{5}{$R[-7]$}{$R[-6]\oplus R[-1]$}
\myrow{$2Z$}{5}{$R[-13]\oplus R[-8]$}{$R[-12]\oplus R[-7]\oplus R[-2]$}
\vskip-5pt
\hbox to\hsize{\hrulefill}
}

\end{document}